\documentstyle[12pt]{amsart} 

\newtheorem{theorem}{Theorem}
\newtheorem{lemma}{Lemma}
\newtheorem{proposition}{Proposition}
\newtheorem{corollary}{Corollary}
\newtheorem{sublemma}{Sublemma}
\newtheorem{remark}{Remark}
\newtheorem{conjecture}{Conjecture}
\newtheorem{definition}{Definition}

\begin{document}
\title{Generalized Bergmann metrics \\ and \\ Invariance of Plurigenera}
\date{February 14, 1996}
\author{Hajime Tsuji}
\address{Hajime Tsuji\\ 
Department of Mathematics\\
Tokyo Institute of Technology\\
2-12-1 Ohokayama Megro 152\\
Tokyo, Japan}
\thanks{Research at MSRI is supported in part by NSF grant DMS-9022140.}
\curraddr{1000 Centennial Drive, Berkeley, CA 94720, U.S.A.} 
\maketitle
\begin{abstract}An invariant kernel for the pluricanonical system  of 
a projective manifold of general type is introduced. 
Using this kernel we prove that the  Yau volume form on a smooth projective variety has seminegative Ricci 
curvature. 
As a biproduct we prove  the invariance of plurigenera for smooth
projective deformations of manifolds of general type. 
\end{abstract}
\tableofcontents
\section{Introduction}

Let $\Omega$ be a bounded domain in ${\bf C}^{n}$ and let 
$\{ \phi_{i}\}_{i=1}^{\infty}$ be a complete orthonormal basis of 
of $L^{2}$ holomorphi $n$-forms on $\Omega$.
Then the Bergmann kernel of $\Omega$ is defined by 
\[
K(z,w) = \sum_{i=1}^{\infty}\phi_{i}(z)\bar{\phi}_{i}(w)
\]
and the Bergmann K\"{a}hler form is defined by 
\[
\omega  = \sqrt{-1}\partial\bar{\partial}\log K(z,z).
\]
The K\"{a}hler metric associated with $\omega$ is called the Bergmann metric
of $\Omega$.
One of the most important property of the Bergmann metric is the invaiance underthe action of the holomorphic automorphisms of $\Omega$. 

In the same way, one can introduce the Bergmann kernel and the Bergmann metric on a complex manifold 
$X$ 
such that the space of $L^{2}$ canonical forms  $H^{0}_{(2)}(X,{\cal O}_{X}(K_{X}))$ is very ample.
In the case of compact complex manifolds, the class of such manifolds corresponds to the class of  complex manifolds
whose canonical bundles are very ample.   Hence it is a relatively small class. 

In this paper, we generalize the Bergmann metric to wider class of manifolds 
and discuss some applications. 
Roughly speaking, we introduce a Bergmann kernel function for pluricanonical systems.

\begin{theorem} Let $X$ be a projective manifold of general type. 
Then for every sufficiently large integer $m$, there exists a singular
K\"{a}hler form $\omega_{m}$ on $X$  and a reproducing kernel  $K_{m}(z,w)$
of holomorphic $m$-ple canonical forms on $X$ with respect to the inner product
induced from the $-m$-th powe of the 
Yau intrisic pseudovolume form  on $X$ such that  
\begin{enumerate} 
\item  $\omega_{m}$ is invariant under $\mbox{Aut}(X)$,
\item  there exists a singular (degenerate) $m$-ple volume form $K_{m}(z,z)$ on 
$X$ such that  
\[
\omega_{m} = \frac{\sqrt{-1}}{m}\partial\bar{\partial}\log K_{m}(z,z)
\]
holds. 
\end{enumerate}  
\end{theorem} 
\begin{remark} To construct $\omega_{m}$ one can use the inner product defined 
by the $-m$-th power of the Kobayashi pseudovolume form. 
Then we have another singular K\"{a}hler form. 
The reason why we use Yau pseudovolume form is mainly the bimeromophic 
invariance of it.  
This enables us to use birational geometry instead of biholomorphic geometry.
But I do not have any specific example which tells us that these pseudovolume
forms are actually different.
\end{remark}

\begin{theorem} Let $X$ be a projective manifold of general type.  
Let $d\mu$ be the
Yau volume form on $X$. 
Then the Ricci curvature of $d\mu$ is  negative  
in the sense of current. 
In fact 
\[ 
d\mu (z) = \limsup_{m\rightarrow\infty}K_{m}(z,z)^{\frac{1}{m}} (z\in X)
\]
and
\[
-\mbox{Ric}\, d\mu  =   \lim_{m\rightarrow \infty}\omega_{m}
\]
hold.
\end{theorem}

As an application of this theorem, we prove the following theorem. 

\begin{theorem}  Let  $\pi : X\longrightarrow S$ be a smooth 
projective family of projective varieties.       
Assume that every fibre of $\pi$ is of general type. 
Then the plurigenera of the fibres are locally constant on $S$.
\end{theorem} 

Theorem 3 is a partial answer to the following conjecture. 
\begin{conjecture} 
The plurigenera is invariant under smooth projective deformations. 
\end{conjecture} 
There have been several works on this conjecture. 

The conjecture is trivial for the family of curves.  
As for the deformation of surfaces, S. Iitaka(\cite{i}) proved the conjecture using
classification of compact complex surfaces.
In general dimensions an interesting approach was proposed by M. Levine (\cite{L}). 
He computed the obstruction for the extention of the member of 
pluricanonical system  on a fibre and proved that it always vanishes 
by using Hodge theory in the case that the member is smooth.  
On the other hand, N. Nakayama pointed out that if the minimal model program 
has been completed (also in the relative case), the conjecture is a direct 
consequence (\cite{na}).  
By the completion of minimal model program for the case of 3-folds, 
J. Kollar and S. Mori proved Conjecture 1 in the case of 3-folds
(\cite[p.535, Theorem 1.3]{k-m}). 
These approaches require some smoothness for a general member of the pluricanonical systems on the fibres.  

Our approach is completely different from the above ones. 
We use the  Yau pseudovolume form on the fibre to  controle the growth of plurigenera and 
use  the continuity of the pseudovolume form  under the deformation.  
To get the disired result we use the theory of analytic Zariski decompositions 
and $L^{2}$-extention theorem in \cite{o-t}.

The author would like to express his hearty thanks to Professor 
S. Kobayashi  who pointed out the reference \cite{w}
and to Professor 
T. Ohsawa for his frank criticism. 
He also would like to express his hearty thanks to Professor S. Lu who
explained Yau pseudovolume form to him. 
Finally he would like to express his hearty thanks to MSRI for their 
hospitality and the stimulating atomosphere. 
\begin{remark} Recently  T. Mabuchi and I have proved the conjecture in full generality
by a different method (in preparation).  
Hence Theorem 3 has beed generalized. 
But the proof of Theorem 3 would be of independent interest. 
\end{remark} 

\part{Generalized Bergmann Kernels} 

\section{Generalized Bergmann metrics}

\subsection{Yau pseudovolume form}
In \cite{ya}, S.-T. Yau introduced an intrinsic pseudovolume form on a complex manifold
which is very  similar to Kobayashi volume form.

Let $M$ be a $n$-dimensional connected complex manifold.
Let $\Delta^{n}$ denote the unit open polydisk in ${\bf C}^{n}$.
Let us take a point $x$ on $M$.
Let $f : \Delta^{n}\longrightarrow M$ be a meromorphic mapping which satisfies
the following three conditions:
\begin{enumerate}
\item $f$ is holomorphic on a neighbourhood of $O$,
\item $f(O) = x$,
\item $f$ is nondegenerate at $O$.
\end{enumerate}
Then we define a pseudovolume form $d\lambda (x)$ at $x$ by
\[
d\mu (x) = \inf (f^{-1})^{*}d\mu_{\Delta^{n}}
\]
where the infimum is taken with respect to all the $f$ which satisfies the
above three conditions and $d\mu_{\Delta^{n}}$ is the Poincar\'{e} volume form 
on $\Delta^{n}$ defined by 
\[
d\mu_{\Delta^{n}} = (\prod_{i=1}^{n}\frac{4}{(1 -\mid\! z_{i}\!\mid^{2})^{2}}
)(\sqrt{-1})^{n}dz_{1}\wedge d\bar{z}_{1}\wedge\cdots\wedge dz_{n}\wedge d\bar{
z}_{n}
\]
Let $f : X -\cdots\rightarrow Y$ be a meromorphic mapping between complex manifolds. 
Then  we have the inequality:
\[
f^{*}d\mu_{Y} \leq d\mu_{X} 
\]
by the definition of the pseudovolume form. 
We call this property the volume decreasing property of Yau pseudovolume form.
By this property it is clear that Yau pseudovolume form is invariant under
bimeromorphic mappings.

It is well known that $d\mu$ is an uppersemicontinuous pseudovolume
form on $X$.
It is clear that Yau pseudovolume form is smaller than Kobayashi pseudovolume 
form  (I do not know any example which shows these volume forms are actually
different).
We call a complex manifold $X$ to be meromorphically measure hyperbolic
if for any nonempty open subset of $X$ the measure defined by $d\mu$
is nonzero.
It is easy to see that  every Kobayashi hyperbolic manifolds are  meromorphicall
y measure hyperbolic. 

\begin{theorem}(\cite{ya})
Let $X$ be a projective manifold of general type.
Then $d\mu$ is nondegenerate on a nonempty Zariski open subset
of $X$.  
In particular $X$ is meromorphically measure hyperbolic. 
\end{theorem}
\begin{lemma}
Let $X$  be a projective manifold of general type. 
Then there exists a proper subvariety $V$ of $X$  such that 
$d\mu$ is continuous on $X - V$. 
\end{lemma}
{\em Proof}.  This fact is essentially in \cite{ya} and was used in \cite{w}
implicitly (see the proof of \cite[p. 367, Theorem 8]{w}).
The proof goes as follows.

For the first by the definiciton  we see that $d\mu$ is uppersemicontinuous. 
Hence we only need to prove that $d\mu$ is lowersemicontinuous. 
If we take $m$ sufficiently large, we may assume that $\mid mK_{X}\mid$ gives a 
birational rational mapping into a projective space. 
Let $\sigma_{0},\ldots ,\sigma_{N}$ be  global holomorphic sections
of $mK_{X}(m >> 1)$ which give a birational rational mapping  $\Phi$ into
${\bf P}^{N}$ and the pullback of the Fubini-Study form 
\[
\sqrt{-1}\partial\bar{\partial}\log\sum_{i=0}^{N}\mid\sigma_{i}\mid^{2}
\]
is a strictly positive on $X$. 
The existence of  $\sigma_{0},\ldots ,\sigma_{N}$ follows from    
Kodaira's lemma (cf.[Appendix]\cite{k-o}). 
Let $V$ denote the support of the  base locus of  $\sigma_{0},\ldots ,\sigma_{N}$. 
Let $f_{i} : \Delta^{n}\longrightarrow X$ be a sequece of meromorphic mappings
such that
\begin{enumerate}
\item $f_{i}$ is holomorphic near the origin and $\{ f_{i}(O)\}$ converges to
a point $x$ on $X - V$,
\item the Jacobian determinant $\mbox{Jac}\,df_{i,O}$ are not zero and converges to a nonzero element. 
\end{enumerate}
We note that $f_{i}^{*}\sigma_{j} (0\leq j\leq N)$ are holomorphic pluricanonica
l forms on $\Delta^{n}$ by Hartogs extention theorem.
Then by the maximum principle, there exists a constant $C$ independent of $i$
such that
\[
(\frac{(\sqrt{-1})^{m\frac{n(n-1)}{2}}f_{i}^{*}(\sum_{j=0}^{N}\sigma_{j}\wedge\bar{\sigma}_{j})}{d\mu_{\Delta^{n}}^{m}})\leq C
\]
holds. 
Here as usual we apply the maximum principle for 
\[
(\frac{(\sqrt{-1})^{m\frac{n(n-1)}{2}}f_{i,\varepsilon}^{*}(\sum_{j=0}^{N}\sigma_{j}\wedge
\bar{\sigma}_{j})}{d\mu_{\Delta^{n}}^{m}}),
\]
where 
\[
f_{i,\varepsilon}(z) 
= f_{i}((1-\varepsilon)z) \,\,\,\,((1-\varepsilon )z = ((1-\varepsilon )z_{1},\ldots , (1-\varepsilon )z_{n}))  \,\,\,\, (0 < \varepsilon << 1)
\]
and letting $\varepsilon$ tend to $0$. 

Then by the normal family argument (Montel's theorem), we see that 
there exists a subsequence $\{ f_{k}\}$ of $\{ f_{i}\}$ such that 
$\{ f_{k}^{*}\sigma_{j}\} (j= 0,\ldots, N)$ converges to a holomorphic $m$-ple canonical forms on $\Delta^{n}$ compact uniformly. 
Then $\{\Phi\circ f_{k}\}$ converges to a meromorphic mapping 
from $\Delta^{n}$ into ${\bf P}^{N}$. 
By the biratonality of $\Phi$
$\{ f_{k}\}$ of $\{ f_{i}\}$  converges to a
meromorphic mapping $f_{\infty} : \Delta^{n}\longrightarrow X$ such that
\[
df_{\infty ,O} =
\lim_{k\rightarrow\infty}df_{k,O}
\]
holds.
We note that by the assumption $x\in X - V$ and the nondegeneracy condition
of the Jacobian, $f_{\infty}$ is holomorphic at $O$.
Then by this fact it is easy to see that $d\mu$ is continuous on $X - V$.
Q.E.D.

\vspace{10mm}

The following corollary is an easy consequence of the above theorem.

\begin{corollary}
Let $X$ be a projective manifold of general type. 
Then there exits a proper subvariety $V$ in $X$ such that 
for every $x\in X - V$, 
there exists a meromorphic mapping $f : \Delta^{n}\longrightarrow X$
which satisfies the following three conditions:
\begin{enumerate}
\item $f$ is holomorphic on a neighbourhood of $O$,
\item $f(O) = x$ holds,
\item $(f^{-1})^{*}(d\mu_{\Delta^{n}})(x) = d\mu_{X}(x)$ holds.
\end{enumerate}
\end{corollary}

\subsection{Construction of the generalized Bergmann kernels}

In this subsection we shall prove Theorem 1. 
Because in some applications, it is more convinient to construct the generalized
Bergmann kernel on the universal covering space, 
we shall consider not only the projective manifolds of 
general type but also its universal covering spaces. 
The proof we present here is for the universal coverings. 
If we restrict ourselves to the case of the manifold itself, the proof  is much more straightforward.  
But it seems to be natural to consider the universal covering in connection withapplications. 
The following proposition follows from the definition of Yau pseudovolume 
form. 

\begin{proposition} Let $M$ be a complex manifold and
let $\pi : \tilde{M}\longrightarrow M$ be an unramified covering.
Then $d\mu_{\tilde{M}} = \pi^{*}d\mu_{M}$ holds.
In particular $M$ is meromorphically  measure hyperbolic if and only if $\tilde{M}$ is
meromorphically measure hyperbolic.
\end{proposition}
	      
Let $X$ be a meromorphically measure hyperbolic manifold.  
We define the subspace ${\bf H}_{m} = {\bf H}_{m}(X)$ of $H^{0}(X,{\cal O}_{X}(mK_{X}))$ by 
\[
{\bf H}_{m} = \{\eta \in H^{0}(X,{\cal O}_{X}(mK_{X}))\mid \,\,
(\sqrt{-1})^{m\frac{n(n-1)}{2}}\int_{X}\frac{\eta\wedge\bar{\eta}}{d\mu_{X}^{m}}d\mu_{X} <\infty \} .
\]
We introduce an inner product on ${\bf H}_{m}$ by
\[
(\eta ,\tau ) = (\sqrt{-1})^{m\frac{n(n-1)}{2}}\int_{X}\frac{\eta\wedge\bar{\tau}}{d\mu_{X}^{m}}d\mu_{X}.
\]
The inner product exists by Schwarz inequality.

\begin{lemma} ${\bf H}_{m}$ is a Hilbert space.
\end{lemma}
{\em Proof}.  Let $\{\eta_{k}\}$ be a Cauchy sequence in ${\bf H}_{m}$.
Let $f : \Delta^{n}\longrightarrow X$ be an embedding into a relatively compact
subset of $X$. 
We shall prove $\{ f^{*}\eta_{k}\}$ converges.      
By the volume decreasing property of Yau pseudovolume form, we see that $\{ f^{*}\eta_{k}\}$ is a Cauchy seqence 
with respect to the inner product
\[
(\eta ,\tau ) = (\sqrt{-1})^{m\frac{n(n-1)}{2}}\int_{X}\frac{\eta\wedge\bar{\tau}
}{d\mu_{\Delta^{n}}^{m}}d\mu_{\Delta^{n}}. 
\]
Hence $\{ f^{*}\eta_{k}\}$ converges to a $m$-ple canonical form. 
Since $f$ is arbitrary, $\{\eta_{k}\}$ converges. 
Q.E.D.

\vspace{10mm}

Let $\{ \phi_{i}\}$ be a complete orthonormal basis for 
${\bf H}_{m}$. 
We set 
\[
K_{m}(z,w) = \sum_{i}\phi_{i}(z)\bar{\phi}_{i}(w)
\]
and
\[
\omega_{m} = \frac{\sqrt{-1}}{m}\partial\bar{\partial}\log K_{m}(z,z).
\]
$K_{m}(z,w)$ exists since  for every $x\inf X$, $ev_{x} : {\bf H}_{m}\longrightarrow {\bf C}\simeq {\cal O}_{X}(mK_{X})/{\cal M}_{x}$ defined
by
\[
ev_{x}(f)  = f(x)
\]
is a bounded functional. 
If ${\bf H}_{m}$ is very ample, $\omega_{m}$ is a $C^{\infty}$-K\"{a}hler metric
on $X$.  
In this case we call $\omega_{m}$ the $m$-ple Bergmann metric of $X$. 

\vspace{10mm}

Let $X$ be a projective manifold of general type and let $\pi :\tilde{X}\longrightarrow X$ be an universal  covering of $X$.  
Let $\omega_{X}$ be a K\"{a}hler form on $X$. 
Let $H^{0}_{(2)}(\tilde{X},{\cal O}_{\tilde{X}}(mK_{\tilde{X}}))$ be the 
space of $L^{2}$ holomorphic sections of $mK_{\tilde{X}}$ with respect to 
$\pi^{*}\omega_{X}$.  
Then by the standard $L^{2}$-estimates for $\bar{\partial}$ operator, we see 
that for every large $m$, $H^{0}_{(2)}(\tilde{X},{\cal O}_{\tilde{X}}(mK_{\tilde{X}}))$ defines a bimeromorphic mapping on to its image in some infinite
dimensional projective space. 

\begin{lemma} ${\bf H}_{m}(\tilde{X}) = H^{0}_{(2)}(\tilde{X},{\cal O}_{\tilde{X}}(mK_{\tilde{X}}))$ holds.
\end{lemma}
{\em Proof}.  
For the first we shall prove that     
\[
H^{0}(X,{\cal L}^{2}(mK_{X},d\mu_{X}^{-m})) \simeq 
H^{0}(X,{\cal O}_{X}(mK_{X}))
\]
holds for every positive integer $m$, where ${\cal L}^{2}(mK_{X},d\mu_{X}^{-m})$is the sheaf defined by 
\[
{\cal L}^{2}(mK_{X},d\mu_{X}^{-m})(U) = \{
\sigma\in\Gamma (U,{\cal O}_{X}(mK_{X}))\mid 
\,  \frac{\mid\sigma\mid^{2}}{d\mu^{m}}\in L^{1}_{loc}(U)\} .
\]

Let $m$ be a sufficiently large positive  integer such that 
$\mid\! mK_{X}\!\mid$ gives a birational rational mapping from $X$ into 
a projective space. 
Let 
\[
\pi_{m} : X_{m}\longrightarrow X
\]
be a resolution of $\mbox{Bs}\mid\! mK_{X}\!\mid$ and let 
\[
\pi_{m}^{*}\mid\! mK_{X}\!\mid = \mid\! P_{m}\!\mid + F_{m}
\]
be the decomposition into the free part and the fixed part.
Then by the assumption $P_{m}$ is nef and big. 
We identify ${\cal O}_{X_{m}}(P_{m})$ as a subsheaf of 
${\cal O}_{X_{m}}(\pi_{m}^{*}(mK_{X}))$
Then by Kodaira's lemma, we see that there exists an effective ${\bf Q}$-divisor
$E$ such that  $P_{m} -\varepsilon E$ is an ample ${\bf Q}$-divisor 
for every sufficiently small positive rational number $\varepsilon$.
Let $r = r(\varepsilon )$ be a positive integer such that 
$r(P_{m}-\varepsilon E)$ is a very ample Cartier divisor. 
Let $\sigma_{0},\ldots ,\sigma_{N}$ be a basis of  
$H^{0}(X_{m},{\cal O}_{X_{m}}(r(P_{m}-\varepsilon E)))$.
Then for every meromorphic mapping $f :\Delta^{n}\longrightarrow X$ 
we have the inequality
\[
\frac{(\sqrt{-1})^{rm\frac{n(n-1)}{2}}\sum_{j=0}^{N}f^{*}\sigma_{j}\wedge\bar{\sigma}_{j}}{d\mu_{\Delta^{n}}^{rm}} \leq C
\]
holds, where $C$ is a positive constant independent of $f$. 
Hence we have that 
\[
d\mu_{X}\geq \frac{1}{\sqrt[rm]{C}}(\sqrt{-1})^{\frac{n(n-1)}{2}}(\sum_{j=0}^{N}\sigma_{j}\wedge\bar{\sigma_{j}})^{\frac{1}{rm}}
\]
holds on $X$. 
This implies that 
\[
\pi_{m}^{*}{\cal L}^{2}(mrK_{X},d\mu_{X}^{-mr}))
\supseteq 
{\cal O}_{X_{m}}(r(P_{m}-\varepsilon E))
\]
holds. 

Letting $\varepsilon$ tend to $0$, by using the above lower estimate for 
$d\mu_{X}$ we see that 
\[
\pi_{m}^{*}{\cal L}^{2}(mK_{X},d\mu_{X}^{-m})
\supseteq
{\cal O}_{X_{m}}(P_{m})
\]
holds. 
Hence for a sufficiently large $m$ we have the desired isomorphism:
\[
H^{0}(X,{\cal L}^{2}(mK_{X},d\mu_{X}^{-m}))
\simeq
H^{0}(X,{\cal O}_{X}(mK_{X})).
\]
By using the ring structure of the canonical ring  $R(X,K_{X})$,
we conclude that the above isomorphim holds for every positive $m$.

Suppose that there exists an element $\phi$ in $H^{0}_{(2)}(\tilde{X},{\cal O}_{\tilde{X}}(mK_{\tilde{X}}))$  
which is not contained in $\phi\ni {\bf H}_{m}$. 
We set for a positive integer $k > 1$,
\[ 
av(\phi^{k}) = \sum_{\gamma\in\pi_{1}(X)}\gamma^{*}\phi^{k}.
\]
Then as in \cite{t}, $av(\phi^{k})$ exists and defines an element of 
$H^{0}(X,{\cal O}_{X}(mkK_{X}))$.  
Then we have for every $x\in X$
\[
\mbox{mult}_{x}(av(\phi^{k})) \geq \mbox{mult}_{x}\mbox{Bs}\mid\! mkK_{X}\!\mid
\]
holds. 
Let us note that the following simple fact that if  
a sequence of complex numbers $\{ a_{\gamma}\}$ such that  $\sum_{\gamma}a_{\gamma}^{k}$ converges for every integer $k\geq 1$ and 
\[
\sum_{\gamma}a_{\gamma}^{k} = 0 
\]
holds for every $k\geq 1$, then $a_{\gamma} = 0$ for every $\gamma$.  

Hence we have for every $x\in \tilde{X}$
\[
\mbox{mult}_{x}(\phi ) \geq
\lim_{k\rightarrow\infty}k^{-1}\mbox{mult}_{\pi (x)}\mbox{Bs}\mid\! mkK_{X}\!\mid
\]
holds. 
This contradicts the assumption that  $\phi$ is not contained in ${\bf H}_{m}$ by Proposition 1 and the first half of the proof of this lemma.
Q.E.D.

\vspace{10mm}

This lemma implies that the kernel $K_{m}$ exists on $\tilde{X}$. \\
And $\omega_{m} = \sqrt{-1}\partial\bar{\partial}\log K_{m}(z,z)$ defines 
a singular K\"{a}hler form on $\tilde{X}$.
We see that by the  definition that $\omega_{m}$ is invariant under $\mbox{Aut}(\tilde{X})$.  Hence $\omega_{m}$ definies a singular K\"{a}hler form on 
$X$.

\begin{corollary}(\cite{t}) Let $X_{1}$,$X_{2}$ be  compact unramified quotients of a 
complex manifold $\tilde{X}$. 
Then $X_{1}$ is of general type if and only if $X_{2}$ is of general type.
\end{corollary}
{\em Proof}.  Suppose that $X_{2}$ is of general type.      
Then by Theorem 1,  $X_{1}$ carries a singular K\"{a}hler form $\omega_{m}$ 
which is the curvature form of a singular hermitian metric of the canonical
bundle of $X_{1}$. 
Then by the standard $L^{2}$-estimates for $\bar{\partial}$ operator, we see 
that $X_{1}$ is also of general type. 
Q.E.D.

\section{Limit of the Generalized  Bergmann kernels}
In this section we shall denote $K_{m}(x,x)$ by $K_{m}(x)$ for simplicity.
\subsection{First properties of Generalized Bergmann kernels}
\begin{lemma} For every $x\in X$ 
\[
K_{m}(x)= \max_{\parallel\phi\parallel=1}(\phi (z)\bar{\phi}(x))
\]
holds. 
\end{lemma}
{\em Proof} It is clear that 
\[
K_{m}(x) \geq \max_{\parallel\phi\parallel=1}(\phi\bar{\phi})(x)
\]
holds by the definition of $K_{m}$. 
On the other hand, 
\[
K_{m} \leq \max_{\parallel\phi\parallel=1}(\phi\bar{\phi})
\]
follows from Riesz's theorem.  Q.E.D.

\vspace{5mm}

\begin{lemma}
\[
\limsup_{m\rightarrow\infty}\int_{X}K_{m}^{\frac{1}{m}}\leq 
\int_{X}d\mu 
\]
holds.
\end{lemma}
{\em Proof}. 
By H\"{o}lder inequality we see that
\[
\int_{X}K_{m}^{\frac{1}{m}}\leq 
(\int_{X}\frac{K_{m}}{d\mu^{m}}d\mu )^{\frac{1}{m}}\cdot (\int_{X}d\mu)^{\frac{m-1}{m}} 
\]
holds.  
We note that 
\[
\int_{X}\frac{K_{m}}{d\mu^{m}}d\mu = \dim H^{0}(X,{\cal O}_{X}(mK_{X}))
\]
holds. 

By Riemann-Roch theorem, we see that there exists a constant $c$ such that  
\[
\dim H^{0}(X,{\cal O}_{X}(mK_{X})) = cm^{n} + o(m^{n})
\]
holds. 
Hence we have that
\[
\lim_{m\rightarrow\infty}
(\int_{X}\frac{K_{m}(z,z)}{d\mu^{m}}d\mu )^{\frac{1}{m}}
= 1       
\]
holds. 
This implies that 
\[
\limsup_{m\rightarrow\infty}\int_{X}K_{m}(z,z)^{\frac{1}{m}} 
\leq 
\int_{X}d\mu
\]
holds. 
Q.E.D. 
\subsection{Poinwise upper estimate of $d\mu$} 
Let us consider the function
\[
F_{m} = \frac{d\mu^{m}}{K_{m}}
\]
on $X$.   

For the first we shall assume that $K_{X}$ is ample.
In this case there exists a point $x_{0}$ where $F_{m}$ takes its minimum
since $F_{m}$ is continuous. 
We set 
\[
h_{m} = (\frac{1}{F_{m}(z_{0})K_{m}})^{1/m}.
\]
Then $h_{m}$ is a singular hermitian metric on $K_{X}$.

\begin{lemma}(\cite[p.103, Lemma 1.2]{ti})
Let $X$ be a smooth projective manifold and let $(L,h)$ be a positive line
bundle on $X$. Let $x_{0}$ be a point on $X$.  Let $g$ be the K\"{a}hler metric
associated with the curvature of $h$.
Choose a normal coordinate $(z_{1},\ldots ,z_{n})$ at $x_{0}$ such that
$x_{0} = (0,\ldots ,0)$.
We choose a local holomorphic frame $e_{L}$ of $L$ at $x_{0}$ such that
the local representation function $a$ is the hermitian metric $h$ has the
property
\[
a(x_{0}) = 1, da(x_{0}) = 0.
\]
 For an $n$-tuple of integers $(p_{1},\ldots ,p_{n})\in {\bf Z}_{\geq 0}^{n}$
 and  an integer $p^{\prime} > p = p_{1}+\cdots +p_{n}$ be a positive integer
there exists $m_{0} > 0$ such that for $m > m_{0}$, there exists a
holomorphic section $S$ in $H^{0}(X,{\cal O}_{X}(mL))$, satisfying
\begin{enumerate}
\item 
\[
\int_{X}\parallel S\parallel_{h^{m}}dV_{g} = 1,
\]
\item 
\[
S(z) = \lambda_{(p_{1},\cdots ,p_{n})} (z_{1}^{p_{1}}\cdots z_{n}^{p_{n}}+ 
O(\mid z\mid^{2p^{\prime}})e_{L}^{m}(1+ O(\frac{1}{m^{2p^{\prime}}})),
\]
where
\[
\lambda^{-2}_{(p_{1},\ldots ,p_{n})} = \int_{X\ \{\mid\!z\!\mid\leq\frac{\log m}{\sqrt{m}}\}}\mid z_{1}^{p_{1}}\cdots z_{n}^{p_{n}}\mid^{2}a^{m}dV_{g}.
\]
where $dV_{g}$ be the volume form with respect to $g$ and 
$\mid\! z\!\mid = \sqrt{\sum\mid \!z_{i}\!\mid^{2}}$. 
\end{enumerate}
\end{lemma}
We call the above $S$'s  the peak sections of $(L,h)$ at $x_{0}$. 
Let $g$ be the K\"{a}hler metric associated with $\mbox{curv}\,h_{m}$. 
For a positive integer $\ell$ we define the number $p(\ell )$ 
by 
\[
p(\ell ) = \max \{\parallel \sigma\parallel_{h_{m}^{\ell}}\!(x_{0}) \mid
\,\,\sigma\in H^{0}(X,{\cal O}_{X}(\ell K_{X})), 
\int_{X}\parallel\sigma\parallel_{h_{m}^{\ell}}^{2}dV_{g} = 1\} ,
\]
and let $\sigma_{\ell}$ denote the section which attains the maximum. 
Then by Lemma 6, we see that 
\[
\lim_{\ell\rightarrow\infty}p(\ell )^{1/\ell} = 1.
\]
Since $h_{m}\geq d\mu^{-m}$ holds on $X$ by the definition of $h_{m}$, 
by Lemma 6, we see that 
\[
\int_{X-\{\mid\!z\!\mid\leq\frac{\log\ell}{\sqrt{\ell }}\}}
\parallel\sigma_{\ell}\parallel_{d\mu^{-m\ell}}^{2}dV_{g} = O(\frac{1}{\ell^{2}}). 
\]
Since $d\mu$ is continuous and $h_{m}(x_{0}) = d\mu^{-1}(x_{0})$ holds
, we see that 
\[
\lim_{\ell\rightarrow\infty}
(\int_{\{\mid\!z\!\mid\leq\frac{\log\ell}{\sqrt{\ell}}\}}
\parallel\sigma_{\ell}\parallel_{d\mu^{-m\ell}}^{2}dV_{g})^{1/\ell} = 1 
\]
holds. 
By Lemma 4 this implies that    
\[
\limsup_{\ell\rightarrow\infty}K_{\ell}(x_{0})^{\frac{1}{\ell}} \geq d\mu (x_{0})
\]
holds.

On the other hand, for a positive integer $\ell$ we define the number $q(\ell )$
by
\[
q(\ell ) =  \max \{\parallel\tau\parallel_{d\mu^{-\ell}}\!(x_{0}) \mid
\,\,\tau\in H^{0}(X,{\cal O}_{X}(\ell K_{X})), 
\int_{X}\parallel\sigma\parallel_{d\mu^{-\ell}}^{2}d\mu = 1\}
\]
and let $\tau_{\ell}$ denote the section which attains the maximum.
By Lemma 4, we have
\[
q(\ell ) = \frac{K_{\ell}(x_{0})}{d\mu^{\ell}(x_{0})} 
\]
holds. 
By the inequality
\[
d\mu^{-1}\leq h_{m},    
\]
as above 
we see that 
$L^{2}$ norm of $\tau_{\ell}$ with respect 
to $d\mu^{-\ell}$
concentrates around $x_{0}$ as $\ell$ tends to infinity.
In fact let $(z_{1},\ldots ,z_{n})$ be a local coordinate around $x_{0}$ as in Lemma 6 
(with respect to $h_{m}$).  
Let $\varepsilon$ be a sufficiently small positive number.
Let $\rho$ be a nonnegative $C^{\infty}$ function of $ \mid z\mid^{2} = \sum_{i=
1}^{n}\mid z_{i}\mid^{2}$
such that
\begin{enumerate}
\item $\rho\equiv 1$ on $B(O,\varepsilon /2) (= \{ z\in X\mid \,\,\mid z\mid\leq
\varepsilon /2\} )$,
\item $\mid d\rho\mid \leq 3/\varepsilon$,
\item $\rho \equiv 0$ on $X - B(0,\varepsilon )$.
\end{enumerate}
If $\varepsilon$ sufficiently small, $\rho$ is well defined.
Let $a$ be a positive function on $B(x,\varepsilon )$ defined as in Lemma 6.
Then there exists a posive constant $c$ such that
\[
a \leq (1 - c\sum_{i=1}^{n}\mid z_{i}\mid^{2})
\]
holds on $B(O,\varepsilon )$. 
This follows from 
the fact that $\mbox{curv}\,h_{m}$ is positive at $x_{0}$. 
On the other hand since $F_{m}$ takes its minimum at $x_{0}$, we see that 
the Taylor expansion of $\tau_{\ell}$ around $x_{0}$ is of the form 
\[
\tau_{\ell} = (c_{\ell} + O(\mid z\mid^{2})\sigma_{m}^{\ell/m}
\]
holds if $m\mid \ell$, where $c_{\ell}$ is a constant  depending on $\ell$. 
For the case that $m$ does not divide $\ell$ we have the similar expansion
of the hermitian norm of $\tau_{\ell}$ with respect to $h_{m}^{\ell}$. 

We set 
\[
G: = \bar{\partial}(\rho \sigma_{\ell}).
\]
Let $\psi$ be a function on $X$ defined by
\[
\psi = n\rho\log (\sum_{i=1}^{n}\mid z_{i}\mid^{2})
\]
Then by the Taylor expansions of $\tau_{\ell}$ and $a$, we see that
there exist a positive constant $C$ independent of $\ell$  and a positive number
$0 <\beta < 1$ such that
\[
\int_{X} e^{-\psi}\parallel F\parallel_{h_{m}^{\ell}}^{2}dv \leq C\cdot c_{\ell}\cdot\beta^{\ell}. 
\]
holds, where  
$dv$ is the pseudovolume form defined by
\[
dv = K_{m}^{\frac{1}{m}}.
\] 
Then by the standard H\"{o}rmander's $L^{2}$ estimates for $\bar{\partial}$
opearator , we have that
there exists  $u\in C^{\infty}(X,\nu K_{X})$ such that
\[
\bar{\partial}u =  G
\]
and
\[
\int_{X}e^{-\psi}\parallel u\parallel_{h_{m}^{\ell}}^{2}dv
\leq
C_{\ell}\int_{X} e^{-\psi}\parallel F\parallel_{h_{m}^{\ell}}^{2}dv \leq C\cdot C_{\ell}\cdot c_{\ell}\cdot\beta^{\ell}
\]
holds, where $C_{\ell}$ is a positive constant satisfying
\[
C_{\ell} = O(\frac{1}{\ell}).
\]
Moreover since $d\mu^{-1}\leq h_{m}$ we see that 
\[
\int_{X}e^{-\psi}\parallel u\parallel_{d\mu^{-\ell}}^{2}dv
\leq  C\cdot
C_{\ell}\cdot c_{\ell}\cdot\beta^{\ell}
\]
holds. 
By the construction
\[
\tilde{\tau}_{\ell} := \rho\tau_{\ell} - u 
\]
is a holomorphic $\ell$-ple canonical form on $X$ and
\[
\tilde{\tau}_{\ell}(x_{0}) = \tau_{\ell}(x_{0})
\]
holds.
Since $\beta < 1$, we see that 
\[
\lim_{\ell\rightarrow\infty}\frac{\parallel\tau_{\ell}\parallel}{\parallel\tilde{\tau}_{\ell}\parallel} = 1
\]
holds.

Then it is clear that 
\[
\limsup_{\ell\rightarrow\infty}q(\ell)^{1/\ell} \leq  1
\]
holds.  
Hence we have that 
\[
\limsup_{\ell\rightarrow\infty}K_{\ell}(x_{0})^{\frac{1}{\ell}} =  d\mu (x_{0})
\]
holds. 
Since $m$ can be arbitrary large we have that
\[
d\mu \geq \limsup_{\ell\rightarrow\infty}K_{\ell}^{\frac{1}{\ell}}
\]
holds.
This implies  that the $L^{2}$-norm of $\tau_{\ell}$ concentrates around
$x_{0}$ as $\ell$ tends to infinity. 
This means that as for the asymptotics, we can work only on a small 
neighbourhood of $x_{0}$. 
By using this fact, we see that
\[
\limsup_{\ell\rightarrow\infty}q(\ell ) \leq 1 
\]
holds. 

Combining the above arguments we see that
\[
d\mu (x_{0}) = \limsup_{\ell\rightarrow\infty}K_{\ell}(x_{0})^{1/\ell}
\]
holds.  
Letting $m$ tend to infinity, we see that 
\[
d\mu \geq \limsup_{\ell\rightarrow\infty}K_{\ell}^{1/\ell}
\]
holds on $X$.

If $K_{X}$ is not ample, we need to modify the argument as follows. 
By Kodaira's lemma (cf. \cite[Appendix]{k-o}) there exists an effective 
${\bf Q}$-divisor $E$ such that $K_{X}- E$ is ample. 
We may assume that the support of $E$ is $V$ in Lemma 1. 
Let $r$ be a positive intger such that $rE$ is a ${\bf Z}$-divisor. 
For every $m \geq  r$ we define the modified generalized Bergmann kernel 
$\tilde{K}_{m}(z,w)$ by
\[
\tilde{K}_{m}(z,w) = \sum_{i}\tilde{\phi}_{i}(z)\bar{\tilde{\phi}}_{i}(w),
\]
where $\{\tilde{\phi}_{i}\}$ is a orthonormal basis of 
$H^{0}(X,{\cal O}_{X}(mK_{X}- rE))$ with respect to the inner product introduced
in 2.2. 
By Lemma 4 one can 
show easily that for $x\in X - V$ 
\[
\lim_{m\rightarrow\infty}(\frac{\tilde{K}_{m}(x)}{K_{m}(x)})^{1/m} = 1
\]
holds by looking at the product:
\[
H^{0}(X,{\cal O}_{X}(m_{0}K_{X}- rE))\times H^{0}(X,{\cal O}_{X}(mK_{X}))
\rightarrow H^{0}(X,{\cal O}_{X}((m + m_{0})K_{X}- rE)),
\]
where $m_{0}$ is a sufficiently large positive integer. 
It is clear that for a sufficiently large $m$,  $d\mu^{m}/\tilde{K}_{m}$ is infinity on $V$ 
and the Ricci form of $\tilde{K}_{m}$ is strictly negative. 
Using $\tilde{K}_{m}$ instead of $K_{m}$, we can argue exactly the same way.
Letting $m$ tend to infinity, we obtain the following lemma.
\begin{lemma}
\[
d\mu \geq \limsup_{\ell\rightarrow\infty}K_{\ell}^{\frac{1}{\ell}}
\]
holds on $X$.
\end{lemma}

\subsection{Lower estimate of the limit of Bergmann kernel}

If we have the pointwise inequality 
\begin{equation}
\limsup_{m\rightarrow\infty}K_{m}(x,x)^{\frac{1}{m}}
\geq
d\mu (x)
\end{equation}
for every $x\in X$, by Lemma 5 , we have the desired equality 
\[
\limsup_{m\rightarrow\infty}K_{m}(x)^{\frac{1}{m}}
= 
d\mu (x)
\] 
for every $x\in X$.
Hence  we shall prove (1).  
Let $m$ be a sufficiently large positive integer.
Let $V$ be a proper subvariety of $X$ as 
in Lemma 1 and  
let $x$ be a point on $X - V$.   
Let $\parallel \,\,\,\,\parallel_{m}$ denote the hermitian norm 
on $mK_{X}$ with respect 
to the hermitian metric $K_{m}^{-1}$  and let  $\parallel\,\,\,\,\parallel$
denote the hermitian norm with respect to $d\mu^{-m}$. 

Let $\sigma\in H^{0}(X,{\cal O}_{X}(mK_{X}))$ be a section of $mK_{X}$ such 
that $\parallel \sigma\parallel_{m} = 1$ and $(\sigma ,\sigma )= 1$ where
$(\,\, ,\,\, )$ denotes the inner product with respect to the hermitian metric
$d\mu^{-m}$ as in the last section. 
Then by Lemma 4, $x$ is the point where $\parallel\sigma\parallel$ takes 
its maximum. 
Hence we have that 
\[
d\parallel\sigma\parallel (x) = 0 
\]
holds. 
Let $f :\Delta^{n}\longrightarrow X$ be a meromorphic  mapping such that
\begin{enumerate}
\item $f$ is holomorphic on a neighbourhood of $O$,
\item $f(O) = x$,
\item $(f^{-1})^{*}d\mu_{\Delta^{n}} = d\mu (x)$.
\end{enumerate} 
Such $f$ exists by Corollary 1.  
If we restrict $f^{-1}$ on a small neighbourhood of $x$, we may consider 
$f^{-1}$ as a local coordinate around $x$.  
We shall denote the coordinate by $(z_{1},\ldots ,z_{n})$. 
We note that the inequality 
\[
d\mu \leq \prod_{i=1}^{n}\frac{4\sqrt{-1}dz_{i}\wedge d\bar{z}_{i}}{(1 - \mid z_{i}\mid^{2})^{2}}
\]
holds on a neighbourhood  of $x$ by the definition of $d\mu$ and the equality
holds at $x$. 
Let  $g_{m}$ be the hermitian metric on $mK_{X}$ on the neighbourhood 
defined by 
\[
g_{m} = (\prod_{i=1}^{n}\frac{4\sqrt{-1}dz_{i}\wedge d\bar{z}_{i}}{(1 - \mid z_{i}\mid^{2})^{2}})^{-m}. 
\]
Then it is clear that 
\[
d\parallel\sigma\parallel_{g_{m}}(x) = 0
\]
holds. 
This implies that 
\[
d(\frac{\prod_{i=1}^{n}(1-\mid z_{i}\mid^{2})^{2m}\sigma\wedge\bar{\sigma}}{\prod_{i=1}^{n}(dz_{i}\wedge d\bar{z}_{i})^{m}})(x) = 0
\] 
holds. 
Hence if we write $(dz_{1}\wedge\cdots\wedge dz_{n})^{m}$ as 
\[
(dz_{1}\wedge\cdots\wedge dz_{n})^{m} = b(z)\sigma  \,\,\,\,\,\,(f\in{\cal O}_{X,x})
\]
on a neighbourhood of $x$, we have that 
\[
b(z_{1},\ldots ,z_{n}) = c + \sum_{i,j}c_{ij}z_{i}z_{j} + O(\mid z\mid^{3})
\]
holds, where $c,c_{ij}$ are constants. 
Let $\varepsilon$ be a sufficiently small positive number. 
Let $\rho$ be a nonnegative $C^{\infty}$ function of $ \mid z\mid^{2} = \sum_{i=1}^{n}\mid z_{i}\mid^{2}$ 
such that 
\begin{enumerate} 
\item $\rho\equiv 1$ on $B(x,\varepsilon /2) (= \{ z\in X\mid \,\,\mid z\mid\leq \varepsilon /2\} )$, 
\item $\mid d\rho\mid \leq 3/\varepsilon$,
\item $\rho \equiv 0$ on $X - B(0,\varepsilon )$.
\end{enumerate}
If $\varepsilon$ sufficiently small, $\rho$ is well defined. 

Let $a$ be a positive function on $B(x,\varepsilon )$ defined by 
\[
a := \parallel\sigma\parallel_{m}^{2}. 
\]
We note that $\limsup_{m\rightarrow\infty} K_{m}^{-1/m}$ has strictly positive curvature on $X - V$.
This fact follows from for a sufficiently large $\ell$ 
$\tilde{K}_{\ell}^{1/\ell}$ defined in the last subsection has 
a strictly positive curvature on $X$ and satisfies the inequality
\[
\tilde{K}_{\ell}^{1/\ell} \leq C_{\ell}d\mu 
\]
on $X$ for some positive constant $C_{\ell}$ by the result in the last subsection.
Then applying Lemma 6 to the cases $(p_{1},\ldots ,p_{n}) = (1,0,\ldots ,0),\ldots ,(0,\ldots ,0,1)$ we see that 
there exists a posive constant $c$ independent of $m$ such that 
\[
a \leq (1 - c\sum_{i=1}^{n}\mid z_{i}\mid^{2})^{m}
\]
holds for every sufficiently large $m$.   
We set for a positive integer $\nu$ 
\[
F = \bar{\partial}(\rho\frac{(dz_{1}\wedge\ldots\wedge dz_{n})^{\otimes\nu}}{c_{\nu}}),
\]
where $c_{\nu}$ is a positive number defined by 
\[
c_{\nu} =  4^{-n\nu}\int_{\Delta^{n}}\prod_{i=1}^{n}(1- \mid z_{i}\mid^{2})^{\nu-1}d\lambda 
\]
where $d\lambda$ is the usual Lebesgue measure.

Let $\psi$ be a function on $X$ defined by
\[
\psi = n\rho\log (\sum_{i=1}^{n}\mid z_{i}\mid^{2})
\]
Then by the Taylor expansions of $b$ and $a$, we see that  
there exist a positive constant $C$ independent of $\nu$  and a positive number 
$0 <\alpha < 1$ such that 
\[
\int_{X} e^{-\psi}\parallel F\parallel_{m,\nu}^{2}dv \leq C(\frac{d\mu^{m}(x)}{K_{m}(x)})^{\nu/m}\alpha^{\nu}.
\]
holds, where $\parallel\,\,\,\,\parallel_{m,\nu}$ denotes the hermitian norm 
induced by $K_{m}^{-\nu/m}$ and 
$dv$ is the pseudovolume form defined by 
\[
dv = K_{m}^{\frac{1}{m}}.
\]
Then by the standard H\"{o}rmander's $L^{2}$ estimates for $\bar{\partial}$
opearator , we have that 
there exists  $u\in C^{\infty}(X,\nu K_{X})$ such that
\[
\bar{\partial}u = F
\]
and 
\[
\int_{X}e^{-\psi}\parallel u\parallel_{m,\nu}^{2}dv 
\leq 
C_{\nu}\int_{X} e^{-\psi}\parallel F\parallel_{m,\nu}^{2}dv \leq C\cdot C_{\nu}
(\frac{d\mu^{m}(x)}{K_{m}(x)})^{\nu/m}\alpha^{\nu}
\]
holds, where $C_{\nu}$ is a positive constant satisfying
\[
C_{\nu} = O(\frac{1}{\nu}).
\]
By the construction  
\[
S_{m,\nu} := \rho\frac{(dz_{1}\wedge\cdots\wedge dz_{n})^{\otimes\nu}}{c_{\nu}}- u 
\]
is a holomorphic $\nu$-ple canonical form on $X$ and 
\[
S_{m,\nu}(x) = \frac{(dz_{1}\wedge\cdots\wedge dz_{n})^{\otimes\nu}}{c_{\nu}}(x)
\]
holds. 
We note that  
\[
\lim_{\nu\rightarrow\infty}(\frac{(dz_{1}\wedge\cdots\wedge dz_{n})^{\otimes \nu}}{c_{\nu}})^{\frac{1}{\nu}}(O) = d\mu_{\Delta^{n}}(O)
\]
holds. 
Hence if
\[
\frac{d\mu (x)}{\sqrt[m]{K_{m}(x)}}\alpha < 1
\]
holds, by the above estimate and Lemma 7 we see that 
\[
d\mu (x) \leq \limsup_{m\rightarrow\infty}K_{m}(x)^{\frac{1}{m}}
\]
holds.

Let $U$ be the subset of $X - V$ defined by
\[
U : = \{ x\in X - V\mid \, (\limsup_{m\rightarrow\infty}K_{m}^{\frac{1}{m}})(x)\ = d\mu (x)\}.
\]
Then by the proof of  Lemma 7, we see that $U$ is nonempty. 
We note that since $\sqrt{-1}\partial\bar{\partial}\log (\limsup_{m\rightarrow\infty}K_{m}^{\frac{1}{m}})$ is a closed positive current, 
$\limsup_{m\rightarrow\infty}K_{m}^{\frac{1}{m}}$ is uppersemicontinuous on $X$. 
 Then by Lemma 1 we see that $U$ is closed.  
On the other hand by the above consideration we see that $U$ is open. 
Hence we conclude that $U$ is equal to $X- V$. 
This implies that the equality 
\[
\limsup_{m\rightarrow\infty}K_{m}^{\frac{1}{m}} = d\mu
\]
holds on $X - V$. 
Since the curvature of  both $\limsup_{m\rightarrow\infty}K_{m}^{\frac{1}{m}}$ and $d\mu$  express  the first Chern class of $X$. 
We see that 
\[
\limsup_{m\rightarrow\infty}K_{m}^{\frac{1}{m}} = d\mu
\]
holds on $X$.

\section{Invariance of Plurigenera}
In this section we shall prove Theorem 3.  
\subsection{Analytic Zariski decomposition}
The notion of analytic Zariski decompositon (AZD) was introduced by 
the author to study big line bundles on an projective variety. 
It is known such a decomposition exists for any big line bundles on 
projective manifolds (\cite{t1}). 
For the first we shall recall the definition of singular hermitian metrics.
\begin{definition} Let $L$ be a holomorphic line bundle over a complex manifold
$M$.  $h$ is said to be a singular hermitian metric on $L$, if there exists 
a $C^{\infty}$-hermitian metric $h_{0}$ on $L$ and a locally $L^{1}$-function 
$\varphi$ such that 
\[
h = e^{-\varphi}h_{0}
\]
holds. 
\end{definition} 
For a singular hermitian line bundle $(L,h)$ as above , we define the curvature 
$\mbox{curv}\,h$ of $(L,h)$ by 
\[
\mbox{curv}\,h =  \mbox{curv}\,h_{0} + \sqrt{-1}\partial\bar{\partial}\varphi .
\]
We define the sheaf ${\cal L}^{2}({\cal O}_{X}(L),h)$ by 
\[
{\cal L}^{2}({\cal O}_{X}(L),h)(U) = \{\sigma\in \Gamma (U,{\cal O}_{X}(L))\mid
 h(\sigma ,\sigma )\in L^{1}_{loc}(U)\} ,
\]
where $U$ runs open subset of $M$. 
We call ${\cal L}^{2}({\cal O}_{X}(L),h)$ the sheaf of germs of $L^{2}$ holomorphic sections of $(L,h)$.  By a theorem of Nadel \cite{n}, we see that 
 ${\cal L}^{2}({\cal O}_{X}(L),h)$ is coherent if the $\mbox{curv}\,h$ is 
 bounded from below by a minus of some K\"{a}hler form locally. 

We shall define an analogy of Zariski decompositions for line bundles 
over a compact complex manifold. 

\begin{definition}
Let $(L,h)$ be a singular hermitian lins bunndle over a compact complex manfold
$M$. 
We call $(L,h)$  an analytic Zariski decomposition, if it satisfies 
the following conditions.
\begin{enumerate}
\item $\mbox{curv}\, h$ is a closed positive $(1,1)$-current on $M$,
\item  the natural morphism $H^{0}(X,{\cal O}_{X}(mL)) \rightarrow 
H^{0}(X,{\cal L}^{2}({\cal O}_{X}(mL),h^{m}))$ is an isomorphim for 
every positive integer $m$.
\end{enumerate}
\end{definition}
\begin{theorem} Let $X$ be a smooth projective variety of general type
and let $d\mu$ be the Yau pseudovolume form on $X$. 
Then $(K_{X},d\mu^{-1})$ is an analytic Zariski decomoposition of $K_{X}$.
\end{theorem}
{\em Proof}. 
By Theorem 2, $\mbox{curv}\, d\mu^{-1}$ is a closed positive current.
By Lemma 3 (and its proof), we have the isomorphism:
\[
H^{0}(X,{\cal L}^{2}({\cal O}_{X}(mK_{X}),d\mu^{-m})) 
\simeq 
H^{0}(X,{\cal O}_{X}(mK_{X}))
\]
for every positive $m$.
Q.E.D.

\subsection{Growth estimate for the sections}

\begin{definition} Let $L$ be a holomorphic line bundle over a compact 
complex $n$-manifold $M$. 
Let $gw(L)$ be the number defined by
\[
gw(L) = \limsup_{m\rightarrow \infty}m^{-n}\dim H^{0}(M,{\cal O}_{M}(mL)).
\]
We call $gw(L)$ the growth of $L$.
\end{definition}

Let $X$ be a smooth projective variety of general type and let 
$d\mu$ be the Yau pseudovolume form on $X$. 
Let $\Theta_{\mu}$ denote the curvature of $d\mu^{-1}$. 
Then we have the Lebesgue decomposition:
\[
\Theta_{\mu} = \Theta_{\mu,abc} + \Theta_{\mu,sing},
\]
where $ \Theta_{\mu,abc}$ denotes the absolutely continuous part of 
$\Theta_{\mu}$ and $\Theta_{\mu,sing}$ denotes the singular part of 
$\Theta_{\mu}$.  
We shall study the meaning of this decomposition. 

\begin{theorem}
\[
\int_{X}(\Theta_{\mu,abc})^{n}  
= (2\pi )^{n}n!gw(K_{X})
\]
holds. 
\end{theorem}
{\em Proof}.  Let $\pi_{m} : X_{m}\longrightarrow X$ be a resolution of 
the base locus of the complete linear system  $\mid\! mK_{X}\!\mid$. 
Let 
\[
\mid\!\pi_{m}^{*}mK_{X}\!\mid = \mid\! P_{m}\!\mid  + F_{m}
\]
be the decomposition of $\mid\!\pi_{m}^{*}mK_{X}\!\mid$ into the free part and 
the fixed part. 
We shall take $\{\pi_{m!}\}$ such that $\pi_{(m+1)!}$ factors throgh $\pi_{m!}$.
We shall write $\pi_{(m+1)!}$ as 
\[
\pi_{(m+1)!} = \xi_{m+1}\circ\pi_{m!}
\]
Then it is clear that 
\[
P_{(m+1)!}- \xi_{m+1}^{*}(m+1)P_{m!} 
\]
is effective. 
Since $P_{(m+1)!}$  and $P_{m!}$ are nef, we see that 
\[
P_{(m+1)!}^{n} - (\xi_{m+1}^{*}(m+1)P_{m!})^{n} \geq 0
\]
holds for every positive integer $m$. 
This implies that 
\[
\lim_{m\rightarrow\infty}\frac{P_{m!}^{n}}{(m!)^{n}} 
\]
exists.
It is clear that 
\[
gw(K_{X}) \geq \frac{1}{(2\pi)^{n}n!}\lim_{m\rightarrow\infty}\frac{P_{m!}^{n}}{(m!)^{n}}
\]
holds. 
But by the work of Fujita (\cite[p.1,Theorem]{f}) we see that the equality
\[
gw(K_{X}) =  \frac{1}{(2\pi)^{n}n!}\lim_{m\rightarrow\infty}\frac{P_{m!}^{n}}{(m!)^{n}}
\]
holds. 
On the other hand by the definition of $\omega_{m}$ , we have that
\[
\int_{X}\omega_{m,abc}^{n} = m^{-n}P_{m}^{n} 
\]
holds. 
Since   
\[
\Theta_{\mu} = \lim_{m\rightarrow\infty}\omega_{m} 
\]
holds, we have that
\[
\frac{1}{(2\pi )^{n}n!}\int_{X}(\Theta_{\mu,abc})^{n}\leq gw(K_{X}) 
\]
holds. 
On the other hand, the Lelong number  $n(\Theta_{\mu}), n(\omega_{m})$ of $\Theta_{\mu},\omega_{m}$ respectively satisfy
the inequality
\begin{equation}   
n(\Theta_{\mu},x) \leq  n(\omega_{m},x) 
\end{equation}
for every $x\in X$ (see the proof of Lemma 3). 
And we note that the de Rham cohomology class of $\Theta_{\mu}$ and 
$\omega_{m}$ are $2\pi c_{1}(K_{X})$. 
Let $B_{m}$ denote the base locus of $mK_{X}$. 
By the construction the support of $n(\Theta_{\mu},x), n(\omega_{m},x)$ are contained in $B_{m}$. 
We shall compute the integral
\[
I_{m} = \int_{X-B}(\Theta_{\mu}^{n}-\omega_{m}^{n}).
\] 
Then we have 
\[
I_{m} = \int_{X-B}(\sqrt{-1}\partial\bar{\partial}\log\frac{d\mu}{K_{m}^{\frac{1}{m}}})(\Theta_{\mu}^{n-1}+ \Theta_{\mu}^{n-2}\omega_{m}+\cdots 
+ \omega_{m}^{n-1}) 
\]
holds.  
We note that 
\[
\Theta_{\mu}^{n-1}+ \Theta_{\mu}^{n-2}\omega_{m}+\cdots 
+ \omega_{m}^{n-1}
\]
is an absolutely continuous positive form on $X-B_{m}$.
Then by Stokes theorem and (2), we conclued that 
$I_{m}$ is nonnegative. 
Hence we have 
\[
\frac{1}{(2\pi )^{n}n!}\int_{X}\Theta_{\mu,abc}^{n} \geq gw_(K_{X})
\]
holds. 
This completes the proof of Theorem 6. 
Q.E.D. 
\begin{remark} 
The other proof will be given (implicitly), in the proof of Theorem 7 
below (see also the proof of Lemma 11).
\end{remark}

\subsection{Analytic Zariski Decomposition of the family}
Let $\pi : X\longrightarrow S$ be a smooth projective family of general type
variety. 
To prove Theorem 3, we may assume that $S$ is the unit open disk $\Delta$ 
in ${\bf C}$.
For every $t\in S$, we denote the fibre $\pi^{-1}(t)$ by $X_{t}$.
Let $d\mu_{t}$ denote the Yau pseudovolume form on $X_{t}$ and 
let $\Theta_{t}$ denote the curvature of $(K_{X_{t}},d\mu_{t}^{-1})$.

\begin{lemma}There exists a proper subvariety $V$ which does not contain 
any fibre of $\pi$ such that 
$d\mu_{t}$ is a continuous relative pseudovolume form on $X - V$. 
\end{lemma} 
{\em Proof}.  Since $K_{X/S}$ is relativiely big by Kodaira's lemma 
\cite[Appendix]{k-o} , we have a finite
number of global sections $\sigma_{0},\ldots \sigma_{N}$ 
of $R^{0}\pi_{*}{\cal O}_{X}(mK_{X/S})$ 
which generates an  ample subbundle of $mK_{X_{t}}$ for every $t\in\Delta$. 
Then we can apply the same proof as that of  Lemma 1.  This  proves the 
lowersemicontinuity of $d\mu_{t}$. 

The uppersemicontinuity of $d\mu_{t}$ follows from the projectivity of 
$\pi$. 
Since $\pi$ is projective, we may assume that the family $X$ is 
embedded in a projective space ${\bf P}^{N}$ by an very ample 
subsehaf of the relative $m$-ple canonical bundle for some large $m$ as 
above.  Let 
$f_{0} : \Delta^{n} -\cdots \rightarrow X_{0}$  be a meromorphic 
mapping which is holomorphic at $O$ and nondegenerate at $O$. 
We want to construct a meromorphic
extention $f : \Delta^{n}\times \Delta (\rho ) -\cdots \rightarrow X$ 
such that $f(z,0) = f_{0}(z) (z\in \Delta^{n})$ and $f(z,t)\in X_{t}(z\in\Delta^{n})$,  
where $\rho$ is a sufficiently small positive number and $\Delta (\varepsilon )$ is the open polydisk of radius $\rho$.  
We claim that such an extension exists if we replace $f(z)$ by $f((1-\varepsilon )z)$ for every positive number $\varepsilon < 1$.
Let $F_{1},\ldots ,F_{r}\in {\cal O}_{\Delta}(\Delta )[X_{0},\ldots ,X_{N}]$ be a homogeneous polynomials with coefficients in ${\cal O}_{\Delta}(\Delta )$
such that for every $t\in\Delta$ $F_{i,t} := F_{i}(t)(1\leq i\leq r)$ is a 
set of defining equations of $X_{t}$ in ${\bf P}^{N}$,  where we have considered
$F_{i}$ as a holomorphic mapping from $\Delta$ to ${\bf C}[X_{0},\ldots ,X_{N}]$. We set $F := (F_{1},\ldots ,F_{r}) : {\bf C}^{N+1}\times \Delta\longrightarrow {\bf C}^{r}$. 
To construct $f$ is equivalent to find set of holomorphic functions 
$g(z,t) = (g_{0},\ldots ,g_{N}): \Delta^{n}\times \Delta (\varepsilon )
\longrightarrow {\bf C}^{N}$ such that  the equation
\[
(*)\,\,\,\,\,\,\,\, F(g(z,t),t) = 0 
\]
holds with the initial condition
\[
g(z,0) = g_{0}
\]
where $g_{0}$ is a lifting of $f_{0}$ to a holomorphic mapping to ${\bf C}^{N+1}$.
We shall expand the solution $g(z,t)$ as 
\[
g(z,t) = \sum_{\nu =0}^{\infty}g_{\mid\nu}(z)t^{\nu}.
\]
Then $g_{\mid 1}(z)$ satisfies the equation:
\[
\frac{\partial (F_{1},\ldots,F_{r})}{\partial (X_{0},\ldots ,X_{N})}(g_{0})^{t}g_{\mid 1}(z) 
+ \frac{\partial F}{\partial t}(g_{0}) = 0.
\] 
This equation for $g_{\mid 1}(z)$ has a holomorphic solution on $\Delta^{n}$
since the polar set of $g_{0}$ is of codimension $\geq 2$ on $\Delta^{n}$
and the equation $(*)$ is equivalent to a equation with $r=1$ by a 
generic projection to ${\bf P}^{n+1}$. 
The equation for $g_{\mid\nu}$ is of the form:
\[
\frac{\partial (F_{1},\ldots,F_{r})}{\partial (X_{0},\ldots ,X_{N})}(g_{0})^{t}g_{\mid\nu} 
+ \mbox{terms with holomorphic coefficients in $g_{\mid i} (0\leq i\leq \nu-1)$} = 0.
\]

 Inductively we can determine holomorphic functions $g_{\mid\nu} (\nu\geq 1)$.
To assure the convergence of the formal solution, we need to construct 
$g_{\mid\nu} (\nu\geq 1)$ with some estimate. 
Let $\varepsilon$ be a positive number less than 1. 
By using the standard argument (see for example \cite{kod}), 
if we replace $f(z)$ by $f((1-\varepsilon )z)$ we can construct a solution 
$g_{\mid\nu}$  with the estimate
\[
\mid g_{\mid\nu}(z)\mid \leq C\cdot M^{\nu} \,\,\,(z\in \Delta^{n})
\]
for some positive constants $C$ and $M$ independent of $\nu$. 
This assures the convergence of the solution $g(z,t)$ for $\mid t\mid < 1/M$. 

The uppersemicontinuity of $d\mu_{t}$ follows immediately from the existence of 
this extention. 
Q.E.D. 

\vspace{10mm}

\begin{lemma} $\int_{X_{t}}\Theta_{t,abc}^{n}$ is a constant function 
on $t\in S$.
\end{lemma} 
We set 
\[
I(t): = \int_{X_{t}}\Theta_{t,abc}^{n}. 
\]
By  Theorem 6, we have that 
\[
I(t) = gw(K_{X_{t}}) 
\]
holds.  By the well known  uppersemicontinuity theorem of the cohomology (cf. \cite[p.351, Theorem 7.8]{kod2}), we have that
$I(t)$  is constant outside of a countably many points on $S$ and strictly 
larger on these points.

On the other hand by Lebesgue-Fatous' lemma , for every convergent sequence
$\{t_{j}\}$ in $S$, we have that
\[
I(\lim_{\j\rightarrow\infty}t_{j}) \leq \liminf_{\j\rightarrow\infty}I(t_{j})
\]
holds. 
This implies that $I(t)$ is constant on $S$.
Q.E.D.

\subsection{Comparison with the global growth}
Let $\pi : X\longrightarrow S$ be a smooth projective family of 
varieties of general type over the unit open polydisc in ${\bf C}$.
To prove Theorem 3 
we may assume that this family is embedded into an irrational pencil 
\[
\varphi : Y\longrightarrow C       
\]
such that 
\begin{enumerate}
\item $Y$ is smooth,
\item $Y$ is of general type,
\item the ristriction morphism:
\[
H^{0}(Y,{\cal O}_{Y}(mK_{Y}))
\rightarrow 
R_{0}\varpi_{*}{\cal O}_{Y}(mK_{Y})_{t}
\]
is surjective for every $t\in C$.
\end{enumerate}
We may assume the third condition by the semipositivity of the 
direct image of the relative pluricanonical sheaf (\cite[Thorem 1]{ka2}\cite{fu2}).    
Let $d\mu$ be the ristriction of the Yau volume form on $Y$ to $X$  and let $\Theta$ be the curvature
of $(K_{X},d\mu^{-1})$. 
By Theorem 2, we see that $\Theta$ is a closed 
positive current on $X$. 
By the construction we see that  the restriction $\Theta_{\mid t}:\Theta\mid_{X_{t}}$ is well defined for every $t\in S$.

\begin{lemma} For every $t$, we have the equality
\[
I(t) = \int_{X}\Theta_{t,abc}^{n} = \int_{X_{t}}(\Theta_{\mid t})_{abc}^{n}
\]
\end{lemma}
{\em Proof} By Lemma 6, $I(t)$ is a constant function on $S$. 
On the other hand outside a countably many points on $S$, the equality
\[
I(t) =  \int_{X_{t}}(\Theta_{\mid t})_{abc}^{n}
\]
holds by the construction of $\Theta_{\mid t}$ and the upperseimicontinuity
theorem. 
Then by the trivial inequality,
\[
I(t) \geq  \int_{X_{t}}(\Theta_{\mid t})_{abc}^{n}
\]
we completes the proof of Lamma 7.
Q.E.D.

\vspace{10mm}

We set 
\[
{\cal F}_{\mid t}(m) :=  {\cal L}^{2}(mK_{X_{t}},d\mu_{\mid X_{t}}^{-m}).
\]

\begin{lemma}       
\[
{\cal F}_{\mid t}(m) = {\cal L}^{2}(X_{t},d\mu_{t}^{-m})
\]
holds.
\end{lemma} 
{\em Proof}.        
To prove this result we shall use the construction of 
analytic Zariski decompositions via parabolic Monge-Amp\`{e}re equation 
in  Part 2. 
We note that we are considering the case of canonical bundles. 
Then by the construction there exists a K\"{a}hler-Einstein form 
$\omega_{1}$,$\omega_{2}$ on a nonempty Zariski open subset $U$ of $X_{t}$ 
such that
\[
{\cal L}^{2}(X_{t},d\mu_{t}^{-m}) =    {\cal L}^{2}(X_{t},(\omega_{1}^{n})^{-m})
\]
and
\[
{\cal F}_{\mid t}(m) =  {\cal L}^{2}(X_{t},(\omega_{1}^{n})^{-m})
\]
holds for every $m \geq 1$. 
Since $\omega_{1},\omega_{2}$ are K\"{a}hler-Einstein on $U$  by the construction of these currents 
 (see Lemma 34,35 in Part 2) we have the pointwise inequality,
\[
\omega_{1}^{n} \geq \omega_{2}^{n}
\]
holds on $U$.
On the other hand by the proof of Theorem 6, we see that 
\[
\int_{U}\omega_{1}^{n} = \int_{U}\Theta_{t,abc}^{n}        
\]
and 
\[
\int_{U}\omega_{2}^{n} = \int_{U}\Theta_{\mid t,abc}^{n} 
\]
hold.
Hence by Lemma 10, we have the equality,
\[
\int_{U}\omega_{1}^{n} = \int_{U}\omega_{2}^{n} 
\]
holds. 
Combining the pointwise inequality between $\omega_{1}^{n}$ and $\omega_{2}^{n}$,
we have that
\[
\omega_{1} = \omega_{2}
\]
holds. 
Since $(K_{X},(\omega_{i}^{n})^{-1}) (i=1,2)$ are AZD, we complete
the proof of Lemma 11. 
Q.E.D. 

\begin{remark} 
Lemma 11 is the only part which requires the use of singular 
K\"{a}hler-Einstein metrics. 
If one succeed in avoiding this very analytic tool, the proof of 
Theorem 3 will become much easier. 
\end{remark}

\subsection{Use of $L^{2}$-extention theorem}
\begin{lemma}
The restriction morphism
\[
H^{0}(X,{\cal O}_{X}(K_{X}\otimes{\cal L}^{2}({\cal O}_{X}(mK_{X}),d\mu^{-m}))
\rightarrow
H^{0}(X_{t},{\cal O}_{X}(K_{X_{t}}\otimes{\cal F}_{\mid t}(m)))
\]
is surjective for every $t\in S$.
\end{lemma} 
{\em Proof}. 
This directly follows from the $L^{2}$-extention theorem 
of \cite{o-t}.  
By Kodaira's lemma, we see that there exists a effective ${\bf Q}$-divisor 
on $X$, such that $K_{X}- E$ is relatively ample on $X$. 
Then there exists an integer $r$ such that $r(K_{X} - E)$ is a very ample 
Cartier divisor on $X$.   
Let $H$ be a smooth member of $\mid r(K_{X} - E)\mid$. 
Then $X - H - \mbox{Supp}\,E$ is Stein and moreover $K_{X}$ is trivial on 
$X - H -  \mbox{Supp}\,E$.  
Now we can apply the extension theorem of \cite{o-t} and completes 
the proof of Lemma 12.
Q.E.D.

\vspace{10mm} 

Now we shall completes the proof of Theorem 3. 
By Theorem 5 and Lemma 11, we see that 
\[
H^{0}(X_{t},{\cal O}_{X}(K_{X_{t}}\otimes{\cal F}_{\mid t}(m))
\simeq
H^{0}(X_{t},{\cal O}_{X_{t}}((m+1)K_{X_{t}}))
\]
holds.                
By Lemma 12, we completes the proof of Theorem 3.

\part{Analytic construction of AZD} 

This part is more or less independent of the last part. 
We shall prove the following theorem.
More general results are in \cite{t1}.      

\begin{theorem} Let $X$ be a smooth projective manifold of general type.
Then there exists a nonempty Zariski open subset $U$ of $X$ and a 
closed positive $(1,1)$-current $\omega_{E}$ such that 
\begin{enumerate} 
\item $\omega_{E}$ is a K\"{a}hler-Einstein form on $U$ with constant 
negative Ricci curvature,
\item $(K_{X},(\omega_{E}^{n})^{-1})$ is an analytic Zariski decomposition.
\end{enumerate}
\end{theorem}

\section{Deformation of K\"{a}hler form I}

In this section we shall consider Hamilton's equation on 
a smooth projective variety of general type and determine 
the maximal existence time for the smooth solution.   
This section is more or less independent from  the other sections.  
Hence a reader who is familiar with Monge-Amp\`{e}re equation may  
skip this section. 
But the basic method of the estimate of parabolic Monge-Amp\`{e}re
 equations is introduced in this section.

\subsection{Hamilton's equation}

Let $X$ be a smooth projective $n$-fold of general type.   
We consider the initial value problem:
\begin{eqnarray}
\frac{\partial\omega}{\partial t}&  = &-\mbox{Ric}_{\omega}- \omega 
 \,\,\,\, 
\mbox{on $X \times [0,T)$} \\
\omega & = &\omega_{0} \,\,\,\, \mbox{on $X\times \{ 0\}$},
\end{eqnarray}
where  \\
$\omega_{0}$ : a $C^{\infty}$-K\"{a}hler form on $X$,\\
$T$ : the maximal existence time for $C^{\infty}$-solution.

Since 
\[
\frac{\partial}{\partial t}(d\omega ) = -d\omega \,\,\, 
\mbox{on $X \times [0,T)$}
\]
\[
d\omega_{0} = 0 \,\,\,\, \mbox{on $X \times \{ 0\}$},
\]
we have that $d\omega = 0$ on $X \times [0,T)$, i.e.,
the equation preserves the K\"{a}hler condition.

\subsection{Reduction to the parabolic Monge-Amp\`{e}re equation} 

Let $\omega$ denote the de Rham cohomology class of
 $\omega$ in $H^{2}_{DR}(X,{\bf R})$. 
Since $-(2\pi )^{-1}\mbox{Ric}_{\omega}$ is a first Chern form of $K_{X}$,
we have
\begin{equation}
[\omega ] = (1- \exp (-t))2\pi c_{1}(L) + \exp (-t)[\omega_{0}].
\end{equation}
Let $\Omega$ be a $C^{\infty}$-volume form on $X$ and let
\[
\omega_{\infty} = curv\,\Omega .
\]
We set
\begin{equation}
\omega_{t} = (1-\exp (-t))\omega_{\infty} + \exp (-t)\omega_{0}.
\end{equation}
Since $[\omega ] = [\omega_{t}]$ on $X\times \{ t\}$ for every $t\in [0,T)$, 
there exists a $C^{\infty}$-function $u$ on $X\times [0,T)$ such that 
\begin{equation}  
\omega = \omega_{t} + \sqrt{-1}\partial\bar{\partial}u. 
\end{equation}
By (3), we have 
\[
\frac{\partial}{\partial t}(\omega_{t} + \sqrt{-1}\partial\bar{\partial}u) 
  =  \sqrt{-1}\partial\bar{\partial}\log 
(\omega_{t}+\sqrt{-1}\partial\bar{\partial}u)^{n}
-(\omega_{t}+\sqrt{-1}\partial\bar{\partial}u) 
\]
Hence
\begin{eqnarray*}
\exp (-t)(\omega_{\infty}-\omega_{0})
+\sqrt{-1}\partial\bar{\partial}(\frac{\partial u}{\partial t}) \\
= \sqrt{-1}\partial\bar{\partial}\log 
(\omega_{t}+\sqrt{-1}\partial\bar{\partial}u)^{n}
- \sqrt{-1}\partial\bar{\partial}\log\,\Omega 
+ \exp (-t)(\omega_{\infty}-\omega_{0}).
\end{eqnarray*}
Then (3) is equivalent to the initial value problem:
\begin{equation}
\frac{\partial u}{\partial t}  =  
\log\frac{(\omega_{t}+\sqrt{-1}\partial\bar{\partial} u)^{n}}{\Omega} - u 
\,\,\,\,\mbox{on $X\times [0,T)$,} 
\end{equation}
where
\[
u  =  0 \,\,\,\,\,\,\,\, \mbox{on $X\times \{ 0\}$}. 
\]

Let
\[
A(X) = \{ [\eta ]\mid \eta :\mbox{K\"{a}hler form on $X$}\}
 \subset H^{2}_{DR}(X,{\bf R})
\]
be the K\"{a}hler cone of $X$.
Since $[\omega ]$ moves on the segument connecting $[\omega_{0}]$ and
 $[\omega_{\infty}] = 2\pi c_{1}(K_{X})$, we cannot expect $T$ to be
 $\infty$, unless $2\pi c_{1}(K_{X})$ is on the closure of $A(X)$ 
in $H^{2}_{DR}(X,{\bf R})$.
We shall determine $T$. 
It is standard to see that $T > 0$ (\cite{ha}).

\begin{theorem}
If $\omega_{0}- \omega_{\infty}$ is a K\"{a}hler form, 
then $T$ is equal to
\[
T_{0} = \sup \{ t>0\mid [\omega_{t}]\in A(X)\} .
\]
\end{theorem}

\vspace{10mm}

The proof of Theorem 8 is almost parallel to that of
 \cite[p.126, Theorem 3]{t}.

\subsection{$C^{0}$-estimate}

\begin{lemma}
If $\omega_{0}-\omega_{\infty}$ is a K\"{a}hler form, then there exists a 
constant $C_{0}$ such that 
\[
\frac{\partial u}{\partial t} \leq C_{0}\exp (-t).
\]
\end{lemma}
{\em Proof}. 
\[
\frac{\partial}{\partial t}(\frac{\partial u}{\partial t}) 
= \Delta_{\omega}\frac{\partial u}{\partial t} 
- \frac{\partial u}{\partial t} 
- \exp (-t)tr_{\omega}(\omega_{0}-\omega_{\infty})
\]
holds by defferentiating (5) by $t$.
By the maximum principle, we have
\[
\frac{\partial u}{\partial t} \leq 
(\max\,\log\frac{\omega_{0}^{n}}{\Omega})\exp (-t).
\]
Q.E.D.

\vspace{10mm}

To estimate $u$ from below, we modify (10) as
\begin{equation}
\frac{\partial u}{\partial t}   =    
\log\frac{(\omega_{t}+\sqrt{-1}\partial\bar{\partial}u)^{n}}{\omega_{t}^{n}} 
+ f_{t}-u \,\,\,\,\mbox{on $X\times [0,T_{1})$,} 
\end{equation}
where 
\[
u   =     0 \,\,\,\,\,\, \mbox{on $X\times \{ 0\}$},
\]
where
\begin{equation}
f_{t} = \log \frac{\omega_{t}^{n}}{\Omega}
\end{equation}
and 
\begin{equation}
T_{1} = \min \{\sup \{ t>0\mid \omega_{t}> 0\} ,T\} .
\end{equation}
If $t\in [0,T_{1})$, we have 
\[
\log \frac{(\omega_{t}+\sqrt{-1}\partial\bar{\partial}u)^{n}}{\omega_{t}^{n}} 
 =  \int_{0}^{1}\frac{d}{ds}\log 
\frac{(\omega_{t}+\sqrt{-1}s\partial\bar{\partial}u)^{n}}{\omega_{t}^{n}}ds \\
  = \int_{0}^{1} \Delta_{s}u ds,
\]
where  $\Delta_{s}$ is the Laplacian with respect to the K\"{a}hler form 
$\omega_{t}+ \sqrt{-1}s\partial\bar{\partial}u$. 
Then by the minimum principle, (11) and Lemma 13, we have
\begin{lemma}
\[
u\geq - C_{0}\exp (-t) + \min_{X}f_{t} \,\,\,
\mbox{on $X\times \{ t\} ,t\in [0,T_{1})$.}
\]
\end{lemma}

\vspace{10mm} 

We note that this estimate is depending on $t$ and $C_{0}$ is 
independent of the choice of $\Omega$. 

\subsection{$C^{2}$-estimate}

For the next we shall obtain a $C^{2}$-estimate of $u$.
\begin{lemma}(\cite[p.351, (2.22)]{y})
Let $M$ be a compact K\"{a}hler manifold and let $\omega$, 
$\tilde{\omega}$ be K\"{a}hler forms on $M$.  
Assume that there exists a $C^{\infty}$-function $\varphi$ such that
\[
\tilde{\omega} = \omega + \sqrt{-1}\partial\bar{\partial}\varphi.
\]
We set 
\[
f = \log \frac{\tilde{\omega}^{n}}{\omega^{n}},
\]
Then for every positive constant such that
\[
C + \inf_{i\neq j}R_{i\bar{i}j\bar{j}} > 1,
\]
\[
\exp (C\varphi )\tilde{\Delta}(\exp (-C\varphi )(n + \Delta\varphi)) \geq 
\]
\[
(\Delta f - n^{2}\inf_{i\neq j}R_{i\bar{i}j\bar{j}}) 
- Cn(n + \Delta\varphi ) 
\]
\[
+ (C + \inf_{i\neq j}R_{i\bar{i}j\bar{j}})
(n + \Delta\varphi )^{\frac{n}{n-1}}\exp (-\frac{f}{n - 1})
\]
holds, where \\
$R_{i\bar{i}j\bar{j}}$ : the bisectional curvature of $\omega$, \\
$\Delta$ : the Laplacian with respect to $\omega$. 
\end{lemma}

Applying this lemma to $\omega_{t}$ and 
$\omega = \omega_{t}+\sqrt{-1}\partial\bar{\partial}u$, we have:
  
\begin{lemma}
For every $C > 0$ depending only on $t\in [0,T_{1})$ such that
\[
C + \inf_{i\neq j}R_{i\bar{i}j\bar{j}}(t) > 1 \;\;\;
\mbox{on $X\times \{ 0\}$},
\]
\[
\exp (Cu)(\Delta_{\omega}-\frac{\partial}{\partial t})
(\exp (-Cu)tr_{\omega_{t}}\omega ) \geq 
\]
\[
-(\Delta_{t}\log \frac{\omega_{t}^{n}}{\Omega} + 
n^{2}\inf_{i\neq j}R_{i\bar{i}j\bar{j}}(t)+n) 
\]
\[
-C(n- \frac{1}{C} -\frac{\partial u}{\partial t})
tr_{\omega_{t}}\omega \\
-\exp(-t)tr_{\omega_{t}}((\omega_{0}-\omega_{\infty})
\cdot\sqrt{-1}\partial\bar{\partial}u) 
\]
\[
+ (C + \inf_{i\neq j}R_{i\bar{i}j\bar{j}}(t))
\exp (\frac{1}{n - 1}(-\frac{\partial u}{\partial t}- u + 
\log 
\frac{\omega_{t}^{n}}{\Omega}))(tr_{\omega_{t}}\omega )^{\frac{n}{n-1}}
\]
holds, on $X\times \{ t\})$ $(t\in [0,T_{1})$, where \\
$\Delta_{t}$ : Laplacian with respect to $\omega_{t}$,\\
$R_{i\bar{i}j\bar{j}}(t)$ : the bisectional curvature  of $\omega_{t}$ \\ 
and $tr_{\omega_{t}}((\omega_{0}-\omega_{\infty})\cdot
\sqrt{-1}\partial\bar{\partial}u)$ denotes the trace
 (with respect to $\omega_{t}$) of the product of the endomorphisms  
$A,B\in \mbox{End}(TX)$ defined by
\[\begin{array}{lll}
\omega_{t}(A(Z_{1})\wedge\bar{Z}_{2}) & = &
 (\omega_{0}-\omega_{\infty})(Z_{1}\wedge\bar{Z}_{2}) \\
\omega_{t}(B(Z_{1})\wedge\bar{Z}_{2}) & = & 
(\sqrt{-1}\partial\bar{\partial}u)(Z_{1}\wedge\bar{Z}_{2}),
\end{array}
\]
where the pair $(Z_{1},Z_{2})$ runs in  $TX\times_{X}TX$. 
\end{lemma}
{\em Proof}. 
Let 
\[
f = \log \frac{\omega^{n}}{\omega_{t}^{n}} \\
\,\,= \frac{\partial u}{\partial t} + u - 
\log \frac{\omega_{t}^{n}}{\Omega}.
\]
Then by Lemma 15, we have
\begin{eqnarray*}
\exp (Cu)\Delta_{\omega}(\exp (-Cu)tr_{\omega_{t}}\omega ) \\
\geq (\Delta_{t}f - n^{2}\inf_{i\neq j}R_{i\bar{i}j\bar{j}}(t))
-Cn(n+\Delta_{t}u) \\ 
+ (C + \inf_{i\neq j}R_{i\bar{i}j\bar{j}})
(tr_{\omega_{t}}\omega )^{\frac{n}{n-1}}\exp (-\frac{f}{n - 1}).
\end{eqnarray*}
Since 
\[
\Delta_{t}f = \Delta_{t}(\frac{\partial u}{\partial t} + u 
- \log \frac{\omega_{t}^{n}}{\Omega}) 
\]
\[
\,\,\,\,\, = \Delta_{t}\frac{\partial u}{\partial t} + 
tr_{\omega_{t}}\omega - n - \Delta_{t}\log \frac{\omega_{t}^{n}}{\Omega}
\]
and 
\[
\exp (Cu)\frac{\partial}{\partial t}(\exp (-Cu)tr_{\omega_{t}}\omega ) 
\]
\[
\,\,\, = -C\frac{\partial u}{\partial t}tr_{\omega_{t}}\omega + 
tr_{\omega_{t}}\frac{\partial \omega}{\partial t} - 
tr_{\omega_{t}}\frac{\partial\omega_{t}}{\partial t}\cdot \omega 
\]
\[
\,\,\, = - C\frac{\partial u}{\partial t}
tr_{\omega_{t}}\omega + \Delta_{t}\frac{\partial u}{\partial t} 
- \exp (-t)tr_{\omega_{t}}(\omega_{0}-\omega_{\infty})
+ \exp(-t)tr_{\omega_{t}}(\omega_{0}-\omega_{\infty})\cdot \omega ,
\]
we obtain the lemma.   Q.E.D. 

\vspace{10mm} 

Let $\varepsilon$ be an arbitrary small positive number.  
We set 
\[
T_{1}(\varepsilon ) =
 \min \{\sup \{ t>0\mid \omega_{t}>0\} -\varepsilon ,T\} 
\]
and let $C$ be a positive number such that 
\[
C + \inf_{i\neq j}R_{i\bar{i}j\bar{j}}(t) > 1 
\]
for all $t\in [0,T_{1}(\varepsilon )]$.
Then since the function $x\exp (-x)$ is bounded on $[0,\infty )$,        
by the maximum principle and Lemma 16, we have that
if $\exp (-Cu)tr_{\omega_{t}}\omega$ take its maximum at 
$(x_{0},t_{0})\in X\times [0,T_{0}(\varepsilon)]$, we have
\[
tr_{\omega_{t}}\omega (x_{0},t_{0}) < C_{\varepsilon}
\]
for some $C_{\varepsilon} > 0$ depending only on $\varepsilon$.
Then by the $C^{0}$-estimate of $u$ , Lemma  13 and Lemma 14, 
by the maximum principle for parbolic equations we have that 
there exists a positive constant $C^{\prime}_{1,\varepsilon}$ such that
\[
tr_{\omega_{t}}\omega < C^{\prime}_{1,\varepsilon}.
\]
Hence we obtain:
\begin{lemma}
There exists a positive constant $C_{1,\varepsilon}$ depending only 
on $T_{1}(\varepsilon )$ such that
\[
\parallel u\parallel_{C^{2}(X)}\leq C_{2,\varepsilon }
\]
for every 
$t\in [0,T_{1}(\varepsilon ))$, where 
$\parallel\,\,\,\,\parallel_{C^{r}(X)}$ is the $C^{2}$-norm
 with repsect to $\omega_{0}$.
\end{lemma}
Now by \cite{tr}, for every $r\geq 2$ there exists 
a positive constant $C_{r,\varepsilon}$ depending only 
on $T_{1}(\varepsilon )$ such that
\[
\parallel u\parallel_{C^{r}(X)} \leq C_{r,\varepsilon}.
\]
Letting $\varepsilon$ tend to 0, we have that
\[
T \geq T_{1}
\]
holds.
Since $[\omega_{T_{0}}]$ is on the closure of the K\"{a}hler cone
 $A(X)$,   
by changing $\Omega$ properly, we can make $T_{0}- T_{1}> 0$
 arbitarary small.
Hence we conclude that $T = T_{0}$.
This completes the proof of Theorem 9.  Q.E.D.

\section{Deformation of K\"{a}hler form II}

In this section we shall construct a K\"{a}hler form on a 
Zariski open subset of X by using a initial value problem similar
to $(3)$ in the last section.
In this section we use the same notation as in the last section.

\subsection{The current $\omega_{E}$}

To state our theorem we need the following definitions.
\begin{definition} Let $D$ be a ${\bf R}$-Cartier divisor on a 
projective variety $Y$. Then the stable base locus of $D$ is defined
by
\[
\mbox{SBs}(D) = \cap_{\nu > 0}\mbox{Supp} \mbox{Bs}\mid [\nu D]\mid .
\]
\end{definition}
\begin{definition} Let $D$ be a Cartier divisor on a projective 
variety $Y$ and let 
$\Phi_{\mid D\mid}; Y -\cdots \rightarrow {\bf P}^{N(\nu )}$ be the 
rational map associated with $\mid\nu D\mid$. 
Let $\mu_{\nu} : Y_{\nu}\longrightarrow Y$ be a resolution of 
the base locus of $\mid\nu D\mid$ and let
 $\tilde{\Phi}_{\mid\nu D\mid};\tilde{X}\longrightarrow 
{\bf P}^{N(\nu )}$ be the associated morphism.
We set
\[
E(\nu D) = 
\overline{\mu_{\nu}(\tilde{E}(\nu D)\cap(Y-\mbox{Supp Bs}\mid\nu D\mid )} 
\,\,\,\mbox{(Zariski closure)}
\]
and call it the exceptional locus of $\mid\nu D\mid$.
It is easy to see that $E(\nu D)$ is independent of the choice of
the resolution of the base locus $\mu_{\nu}$. 
We set
\[
\mbox{SE}(D) = \cap_{\nu >0}E(\nu D)
\]
and call it the stable exceptional locus of $D$.
\end{definition}

\vspace{10mm}
We set
\[
S = \mbox{SBs}(L)\cup SE(L).
\]

\vspace{10mm}

The main result in this section is the following theorem.

\begin{theorem}
There exists a $d$-closed positive $(1,1)$-current $\omega_{E}$ 
on $X$ such that 
\begin{enumerate}
\item $\omega_{E}$ is smooth on a nonempty Zariski open subset $U$ of $X$,  
\item $\omega_{E} =  - \mbox{Ric}_{\omega_{E}}$ 
holds on $U$.
\item $[\omega_{E}] = 2\pi c_{1}(L)$.
\end{enumerate}
\end{theorem}

\subsection{Kodaira's Lemma}

Let $\nu$ be a sufficiently large positive integer such that
\begin{enumerate}
\item 
$\mid\nu L\mid$ gives a birational rational map from $X$ into 
a projective space.
\item 
 Supp Bs$\mid \nu L\mid = \mbox{SBs}(L)$.
\end{enumerate}
Let $f_{\nu} : X_{\nu}\rightarrow X$ be a resolution of 
the base locus of $\mid\nu L\mid$ and let  
\[
F_{\nu} = \sum_{i} b^{\nu}_{i}F_{i}^{\nu}
\]
be the fixed part of $\mid f^{*}_{\nu}(\nu L)\mid$.     
We take $f_{\nu}$ so that $F^{\nu}$ is a divisor
 with normal crosssings.      
We set 
\[
\tilde{b}^{\nu}_{i} = b^{\nu}_{i}/\nu .
\]
Let $\sigma_{i}^{\nu}$ be a global holomorphic section of
 ${\cal O}_{X_{\nu}}(F_{i}^{\nu})$ with divisor $F_{i}^{\nu}$.  
Then there exist hermitian metrics $\parallel\,\,\,\,\parallel$
 on ${\cal O}_{X_{\nu}}(F^{\nu}_{i})$'s  such that 
\[
\omega^{\nu}_{\infty} = f_{\nu}^{*}\omega_{\infty} 
+ \sum\sqrt{-1}\tilde{b}_{i}^{\nu}\partial\bar{\partial}
\log\parallel\sigma_{i}^{\nu}\parallel^{2}
\]
is positive on $f_{\nu}^{-1}(X-S)$, if $\nu$ is sufficiently large.
We may assume 
\[
\log\parallel\sigma_{i}^{\nu}\parallel\leq 0
\]
holds for every  $i$. 
We set

By Kodaira's lemma (\cite[Appendix]{k-o} we have the following lemma 
\begin{lemma}There exists an effective ${\bf Q}$-divisor
\[
R_{\nu} = \sum r_{j}^{\nu}R_{j}^{\nu}
\]
on $X_{\nu}$ such that 
\[
f^{*}_{\nu}(K_{X}) - \sum\tilde{b}_{i}^{\nu}F_{i}^{\nu} - R^{\nu}
\]
is an ample ${\bf Q}$-divisor on $X_{\nu}$.
\end{lemma} 
We note that $\varepsilon R_{\nu}$ has the same property as 
$R_{\nu}$ for $\varepsilon\in [0,1]$.
Let $\tau_{j}^{\nu}$ be a global section of
 ${\cal O}_{X_{\nu}}(R_{j}^{\nu})$ with divisor $R_{j}^{\nu}$. 
Then there exists hermitian metrics $\parallel\,\,\,\parallel$ on
${\cal O}_{X_{\nu}}(R_{j}^{\nu})$ such that  
\[
\omega^{\nu}_{\infty} 
+ \sum_{j}\sqrt{-1}r_{j}^{\nu}\partial\bar{\partial}
\log\parallel\tau^{\nu}_{j}\parallel^{2}
\]
is a smooth K\"{a}hler form on $X_{\nu}$ and 
$\parallel\tau^{\nu}_{j}\parallel\leq 1$ holds on $X_{\nu}$ for all $j$.
We set
\[
\delta_{\nu} 
= \sum_{j}\sqrt{-1}r_{j}^{\nu}\log\parallel\tau^{\nu}_{j}\parallel^{2}. 
\]
Then for every  $\varepsilon\in [0,1]$,
\[
\omega_{\infty}^{\nu} +
 \varepsilon\sqrt{-1}\partial\bar{\partial}\delta^{\nu}
\]
is a smooth K\"{a}hler form on $X_{\nu}$.

We set 
\[
\xi_{\nu} 
= \sum_{i}\tilde{b}_{i}^{\nu}\log\parallel\sigma_{i}^{\nu}\parallel^{2}.
\]
\subsection{Construction of a suitable ample divisor}

To construct K\"{a}hler-Einstein current on $X$, we use the 
Dirichlet problem for parabolic Monge-Amp\`{e}re equation. 
Hence we shall construct a  strongly pseudoconvex  convex exhaustion
of a Zariski open subset of $X_{\nu}$ with certain properties. 
We fix sufficiently large $\nu$ hereafter.            
Let 
\[
\Phi : X_{\nu} \longrightarrow {\bf P}^{N}
\]
be a embedding of $X_{\nu}$ into a projective space.
Let 
\[
\pi_{\alpha} ; X_{\nu} \longrightarrow {\bf P}^{n} (\alpha = 1,\ldots m)
\]
be generic projections and we set
\[
W_{\alpha} : \mbox{the ramification divisor of $\pi_{\alpha}$}, \\
H_{\alpha} :=  \pi_{\alpha}^{*}(\{ z_{0} = 0\} ),
\]
where 
$[z_{0}:\ldots :z_{n}]$ be the homogeneous coordinate of ${\bf P}^{n}$. 
For simplicity  we shall denote the support of a divisor 
by the same notation as the one,
if without fear of confusion.
If $m$ is sufficiently large, we may assume the following conditions:
\begin{enumerate}
\item $\cap_{\alpha = 1}^{m}(W_{\alpha}+H_{\alpha}) = \phi$,
\item $D := (F_{\nu} + \sum_{\alpha = 1}^{m}(W_{\alpha}+ H_{\alpha}))_{red}$ 
is an ample divisor with normal crossings.
\item $D$ contains $S \cup R_{\nu}$.
\item $K_{X_{\nu}}+D$ is ample. 
\end{enumerate}
Then  $X_{\nu} - D$ is strongly pseudoconvex 
and the following lemma is  necessary for our purpose. 
\begin{lemma}
There exists a positive strongly plurisubharmonic exhaustion function
 $\varphi$ of $X_{\nu} - D$ such that
 $\omega_{\varphi} =\sqrt{-1}\partial\bar{\partial}\varphi$
is a complete K\"{a}hler form on $X_{\nu}- D$.
\end{lemma}
{\em Proof}.
Let $D = \sum_{k}D_{k}$ be the irreducible decomposition of $D$ 
and let $\lambda_{k}$ be a global holomorphic section of 
${\cal O}_{X_{\nu}}(D_{k})$ with divisor $D_{k}$. 
Then there exist hermitian metrics $\parallel \,\,\,\parallel$'s on 
${\cal O}_{X_{\nu}}(D_{k})$'s such that     
\[
-\sum_{k}\sqrt{-1}\partial\bar{\partial}
\log \parallel\lambda_{k}\parallel^{2} 
\]
is a smooth K\"{a}hler form on $X_{\nu}$.  
We set for a positive number $\iota$
\[
\varphi = -\sum_{k}\log\parallel\lambda_{k}\parallel
 -\iota\log\log\frac{1}{\parallel\lambda_{k}\parallel}.
\]
Then if we choose $\iota$ sufficiently small, then  
\[    
\omega_{\varphi} = \sqrt{-1}\partial\bar{\partial}\varphi
\]
is a complete K\"{a}hler form on $X_{\nu} - D$. 
Clearly by adding a sufficiently large positive number, we can 
make the exhuastion $\varphi$ to be positive on $X_{\nu} - D$. 
Q.E.D.

\begin{remark} As one see in Lemma 24 below, $\omega_{\varphi}$ has
 a bounded Poincar\'{e} growth.
\end{remark}

\subsection{The Dirichlet problem}

We set for $c > 0$, 
\[
K_{c} = \{ x\in X_{\nu}: \varphi (x) \leq c\} .
\]
 It is easy to see that we may assume that 
there exists a positive constant $c_{0}$ such that the boundary 
$\partial K_{c}$ is smooth for every  $c \geq c_{0}$.
We fix such $c$ and set
\[
K : = K_{c}
\]
for simplicity.  

In the estimate of $u$, we try to make the estimate independent of 
$c \geq c_{0}$ for the later use.

Since $X_{\nu}- D$ is canonically biholomorphic to 
a Zariski open subset of $X$,
we may consider $K$ as a compact subset of $X$. 
We consider the following Dirichlet problem for a parabolic
Monge-Amp\`{e}re equation.
\[ \left\{ \begin{array}{llll}
\vspace{5mm}
\frac{\partial u}{\partial t} & = & 
\log\frac{(\omega_{t}+\sqrt{-1}\partial\bar{\partial}u)^{n}}{\Omega} 
- u & \mbox{on $K\times [0,T)$}   \vspace{5mm} \\
u  & = & (1 - e^{-t^{4}})\xi_{\nu} & 
\mbox{on $\partial K\times [0,T)$} \;\;\mbox{(10)}
\vspace{5mm} \\
u  & = & 0 & \mbox{on $K\times\{ 0\}$} ,
\end{array}
\right.
\]
\addtocounter{equation}{+1}
where $\Omega$ is a smooth volume form on $X$,  
\[
\omega_{t} = (1 - e^{-t^{4}})curv\, h  + e^{-t^{4}}\omega_{0}
\]
and
$T$ is a maximal existence time for the smooth solution on
 $\bar{K}$(the closure in the usual topology).
We shall assume that 
\[
\omega_{0}^{n} = \Omega
\]
holds.    
It is easy to find such $\omega_{0}$ and $\Omega$ by using 
the solution of Calabi's conjecture (\cite{y}). 
We set 
\[
\omega_{\infty} = curv\, h.
\]
By multiplying  a common  sufficiently large positive number to 
$\omega_{0}$ and $\Omega$, if neceessary, we may assume that 
\[
\omega_{0} + \mbox{Ric}\,\Omega = \omega_{0}-\omega_{\infty}
> 0
\]
holds.   
Please do not confuse $u$ with the one in the last section. 
We use the same notation for simplicity. 
For the first we shall show
\begin{theorem}
$T$ is infinite and 
\[
u_{\infty} = \lim_{t\rightarrow\infty}u
\]
exists in $C^{\infty}$-topology on $\bar{K}$. 
\end{theorem} 

\subsection{$C^{0}$-estimate}

We note that by the above choice of $\omega_{0}$,$\Omega$ and the 
Dirichlet condition, the Dirichlet problem (10) is compatible up to 
3-rd order on the corner $\partial K\times\{ 0\}$. 
 
Hence by the standard implicit function theorem, 
we see that $T$ is 
positive. 
Suppose $T$ is finite.  
Then if $\lim_{t\rightarrow T}u$ exist on $K$ in $C^{\infty}$-topology,
 then this is a contradiction.              
Because again by the implicit function theorem, we can continue
the solution a little bit more.  
Hence to prove Theorem 10 it is sufficient to obtain an estimate 
of $C^{k}$-norm of $u$ on $K$ which is independent of $t$. 

We begin with $C^{0}$-estimate. 

\begin{lemma}
There exists a constant $C_{0}^{+}$ such that
\[
\frac{\partial u}{\partial t}\leq C_{0}^{+}e^{-t} \,\,\,\,
\mbox{on $K\times [0,T)$}
\]
holds.
\end{lemma}
{\em Proof}
We set
\[
\omega = \omega_{t} + \sqrt{-1}\partial\bar{\partial}u
\]
and
\[
\tilde{\Delta} = tr_{\omega}\sqrt{-1}\partial\bar{\partial}.
\]
As in the proof of Lemma 3.1, we have
\[
\frac{\partial}{\partial t}(\frac{\partial u}{\partial t}) = 
\tilde{\Delta}\frac{\partial u}{\partial t} - 
\frac{\partial u}{\partial t} 
- 4t^{3}e^{-t^{4}}tr_{\omega}(\omega_{0}-\omega_{\infty})
\]
holds on $K\times [0,T)$. 
Since
\[
u = (1 - e^{-t^{4}})\xi_{\nu} \,\,\,\,\mbox{on $\partial K\times [0,T)$}
\]
and $\omega_{0}- \omega_{\infty}$ is a K\"{a}hler form on $X$, by 
maximal principle
\[
\frac{\partial u}{\partial t} \leq C_{0}^{+}e^{-t} \,\,\,
\mbox{on $K\times [0,T)$}
\]
holds for
\[
C_{0}^{+} = \max\{\max_{x\in K}\log\frac{\omega_{0}^{n}}{\Omega}(x),
\max_{x\in\partial K}\xi_{\nu}(x)\}.
\]
Q.E.D.

\vspace{10mm}

We set 
\[
v = u - (1 - e^{-t^{4}})\xi_{\nu},
\]
\[
\Omega_{\nu} = \exp (\xi_{\nu})\Omega ,
\]
and
\[
\hat{\omega}_{t}
 = \omega_{t} + (1 - e^{-t^{4}})\sqrt{-1}\partial\bar{\partial}\xi_{\nu}.
\]
Then $v$ satisfies the equation:
\[ \left\{ \begin{array}{llll} 
\frac{\partial v}{\partial t} & = &
\log\frac{(\hat{\omega}_{t}
+\sqrt{-1}\partial\bar{\partial}v)^{n}}{\Omega_{\nu}} - v & 
\mbox{on $K\times [0,T)$}  \vspace{8mm}\\
v & = & 0 & \mbox{on $\partial K\times [0,T)$}  \;\; 
\mbox{(11)} \vspace{5mm} \\
v & = & 0 &\mbox{on $K\times\{ 0\}$} \end{array} 
\right. \]
\addtocounter{equation}{+1}
We note that $\hat{\omega}_{t}$ is a K\"{a}hler form on 
$(X_{\nu} - D)\times [0,\infty ]$ and 
$\hat{\omega}_{\infty} = \omega^{\nu}_{\infty}$.  
Then since
\[
\log\frac{(\hat{\omega}_{t}+
 \sqrt{-1}\partial\bar{\partial}v)^{n}}{\hat{\omega}_{t}^{n}} 
= \int_{a=0}^{1}\tilde{\Delta}_{a}v ,
\]
where $\tilde{\Delta}_{a}$ is the Laplacian with respect to 
the K\"{a}hler form 
\[
\hat{\omega}_{t} + a\sqrt{-1}\partial\bar{\partial}v,
\]
by maximum principle and Lemma 19, we have that
\[
v\geq 
\min\{ \min_{x\in K}\log\frac{\hat{\omega}_{t}^{n}}{\Omega_{\nu}}(x) 
- C_{0}^{+}e^{-t},0\} \,\,\mbox{on $K\times [0,T)$}
\]
holds.  
Hence we have
\begin{lemma}
\[
u\geq C_{0}^{-} + (1 - e^{-t^{4}})\xi_{\nu} \,\,\,
\mbox{on $K\times [0,T)$},
\]
where
\[
C_{0}^{-} = 
\min\{ \inf_{(x,t)\in K\times [0,\infty )}\log
\frac{\hat{\omega}_{t}^{n}}{\Omega_{\nu}}(x,t), 0\} - C_{0}
\]
\end{lemma}

\vspace{10mm} 

We note that $C_{0}^{-}$ may depend on $K$ because 
$\log (\hat{\omega}^{n}/\Omega_{\nu})$ may not be bounded from 
below on $X_{\nu} - D$. 
To obtain the $C^{0}$-estimate from below which is independent of 
$K$, we shall consider for $\varepsilon \in (0,1]$, 
\[
v_{\varepsilon} = u - (1- e^{-t^{4}})(\xi_{\nu} + \varepsilon\delta_{\nu}).
\] 
Then by the same argument,  we have
\begin{lemma}
Let $\varepsilon \in (0,1]$.  
Then there exists a constant $C_{0}^{-}(\varepsilon )$ which 
is independent of $K$ such that
\[
u\geq C_{0}^{-}(\varepsilon ) +  
(1 - e^{-t^{4}})(\xi_{\nu}+\varepsilon\delta_{\nu}) \,\,
\mbox{on $K\times [0,T)$}
\]
\end{lemma}

\vspace{10mm}

The reason why $C_{0}^{-}(\varepsilon )$ is independent of $K$ 
is simply because 
\[
\hat{\omega}_{\infty}
 + \varepsilon\sqrt{-1}\partial\bar{\partial}\delta_{\nu}
\]
on $X_{\nu} - D$ extends to a smooth K\"{a}hler form on $X_{\nu}$ 
and
\[
\exp(\xi_{\nu}+\varepsilon\delta_{\nu})f_{\nu}^{*}\Omega
\]
is a smooth semipositive $(n,n)$ form on $X_{\nu}$.

\subsection{$C^{1}$-estimate on $\partial K$}

Hereafter we estimate derivatives of $v$ basically  
by using the method in \cite{c}.
But our estimates is a little bit more complicated 
because we are working on a quasi-projective variety
 which cannot admits a global flat 
K\"{a}hler metric. 
We set 
\[
\psi_{c} = b(\varphi - c),
\]
where $b$ is a positive constant. 
We note that 
\[
\psi_{c} = 0 \,\,\,\mbox{on $\partial K$}
\]
 and 
\[
\psi_{c} < 0  \,\,\,\ \mbox{on $K$}.
\]
Then since $\omega_{\varphi}$ is a complete K\"{a}hler form of 
Poincar\'{e} growth, if we take $b$ sufficiently large 
\[
\log\frac{(\hat{\omega}_{t}
+\sqrt{-1}\partial\bar{\partial}\psi_{c})^{n}}{\Omega_{\nu}} 
- \psi_{c} \geq 0  \,\,\mbox{on $K\times [0,\infty)$}
\]
holds.
It is easy to see that we can take $b$ independent of $c$ and $t$.
Then we have
\[
\left\{  \begin{array}{lll}
\frac{\partial (v - \psi_{c})}{\partial t} & \geq
\log\frac{(\hat{\omega}_{t} + \sqrt{-1}\partial\bar{\partial}v)^{n}}
{(\hat{\omega}_{t} + \sqrt{-1}\partial\bar{\partial}\psi_{c})^{n}}
- (v -\psi_{c}) & \mbox{on $K\times [0,T)$} \vspace{5mm} \\
v - \psi_{c} & = 0 & \mbox{on $\partial K\times [0,T)$} \vspace{5mm} \\
v - \psi_{c} & = -\psi_{c} & \mbox{on $K\times \{ 0\}$}
\end{array} \right.
\]
Since
\[
\log\frac{(\hat{\omega}_{t}+\sqrt{-1}\partial\bar{\partial}v)^{n}}
{(\hat{\omega}_{t}+\sqrt{-1}\partial\bar{\partial}\psi_{c})^{n}}
= \int_{0}^{1}\dot{\Delta}_{a}(v - \psi_{c})da, 
\]
where  $\dot{\Delta}_{a}$ is the Laplacian with respect to the 
K\"{a}hler form 
\[
\hat{\omega}_{t} + \sqrt{-1}\partial\bar{\partial}\{ (1-a)\phi_{c}+av\}
\],
by the  maximum principle, we obtain
\[
v \geq \psi_{c}  \,\,\,\,\,\mbox{on $K\times [0,T)$}.
\]
On the other hand, trivially 
\[
\hat{\Delta}_{t}v \geq -n \,\,\,\mbox{on $K\times[0,T)$}
\]
holds, where $\hat{\Delta}_{t}$ is the Laplacian with respect to
the K\"{a}hler form $\hat{\omega}_{t}$.   
Let $h$ be the $C^{\infty}$-function on $\bar{K}\times [0,T)$ such that
\[\left\{ \begin{array}{llll}
\hat{\Delta}_{t}h & = & -n  &\mbox{on $K\times [0,T)$}  \\     
h & = &0  &\mbox{on $\partial K\times [0,T)$.}
\end{array} \right.
\]
Then by the maximum principle, we have
\[
v\leq h   \,\,\,\,\mbox{on $K\times [0,T)$}.
\] 
Hence we have
\begin{equation}
\psi_{c}\leq v\leq h \,\,\,\,\mbox{on $K\times [0,T)$}.
\end{equation}

Now to fix $C^{k}$-norms on $X_{\nu} - D$, we shall construct 
a complete K\"{a}hler-Einstein  form on $X_{\nu}- D$. 

We quote the following theorem.

\begin{theorem}(\cite{k})
Let $\bar{M}$ be a nonsingular projective manifold and let $B$ be
an effective divisor with only simple normal crossings. 
If $K_{\bar{M}}+ B$ is ample, then there exists a unique
 (up to constant multiple) complete
K\"{a}hler-Einstein form on $M = \bar{M} - B$ with negative Ricci 
curvaure.
\end{theorem} 

By the construction of $D$, $D$ is a divisor with simple normal 
crossings and $K_{X_{\nu}}+ D$ is ample. 
Hence by Theorem 11, there exists a complete K\"{a}hler-Einstein 
form $\omega_{D}$ on $X_{\nu}- D$ such that
\[
\omega_{D} = - \mbox{Ric}_{\omega_{D}}.
\]

Then we have
\[
\parallel dv\parallel \leq \max \{ \parallel dh\parallel, 
\parallel d\psi_{c}\parallel\} \,\,\, 
 \mbox{on $\partial K\times [0,T)$},
\]
where $\parallel\,\,\,\,\parallel$ is the pointwise norm
 with respect 
to $\omega_{D}$. 
To make this estimate independent of $K$, we need to use special
properties of $\omega_{D}$.

\begin{definition} Let $V$ be an open set in ${\bf C}^{n}$. 
A holomorphic map from $V$ into a complex manifold $M$ of dimension 
$n$ is called a {\em quasi-coordinate map}
 iff it is of maximal rank everywhere on $V$. 
$(V;\mbox{Euclidean coordinate of ${\bf C}^{n}$})$ is called a 
{\em local quasi-coordinate} of $M$.
\end{definition}

\begin{lemma}(cf. \cite[p.405, Lemma 2 and pp. 406-409]{k} )
There exists a family of local quasi-coordinates  
${\cal V} = \{ (V;v^{1},\ldots,v^{n})\}$ of $X_{\nu}- D$  with 
the following properties. 
\begin{enumerate}
\item $X_{\nu} - D$ is covered by the images of
 $(V;v^{1},\ldots ,v^{n})$'s.
\item The completment of some open neighbourhood of $D$ is covered 
by a finite number of $(V,v^{1},\ldots ,v^{n})$'s which are
 local coordinate in the usual sense.
\item Each $V$, as an open subset of the complex Euclidean space
 ${\bf C}^{n}$, contains a ball of radius $1/2$.
\item There exists positive constants $c_{D}$ 
and ${\cal A}_{k} (k=0,1,2,\ldots )$ independent of $V$'s 
such that at each  
$(V,v^{1},\ldots v^{n})$, the inequalities:
\[\begin{array}{l}
\frac{1}{c_{D}}(\delta_{ij}) < (g^{D}_{i\bar{j}})
 < c_{D}(\delta_{ij}), \vspace{8mm} \\
\mid (\partial^{\mid p\mid+\mid q\mid}/
\partial v^{p}\partial\bar{v}^{q})g^{D}_{i\bar{j}}\mid
 < {\cal A}_{\mid p\mid +\mid q\mid}, 
\mbox{for any multiindices $p$ and $q$}
\end{array}
\]
hold,  where $g^{D}_{i\bar{j}}$ denote the components of 
$\omega_{D}$
with respect to $v^{i}$'s.
\end{enumerate}
\end{lemma}

\begin{definition} $(\bar{M},B)$ be a pair of smooth projective variety 
of dimension $n$ and a divisor with simple normal crossings on it. 
A complete K\"{a}hler metric $\omega_{M}$ on $M = \bar{M} - B$ 
is said to have bounded Poincar\'{e} growth on
$(M,B)$ if for any polydisk 
$\Delta^{n} = \{ (z_{1},\ldots ,z_{n})\in {\bf C}^{n}\mid 
\mid z_{i}\mid < 1 (1\leq i\leq n)\}$ in $\bar{M}$ such that
\[
\Delta^{n} \cap B 
= \{ (z_{1},\ldots ,z_{n})\in \Delta^{n}\mid \,\,\,
 z_{1}\cdots z_{k} = 0\} (k\leq n),
\]
$\omega_{M}\mid\Delta^{n}$ is quasi-isometric to 
\[
\omega_{P} = \sum_{i=1}^{k}
\frac{\sqrt{-1}dz_{i}\wedge d\bar{z}_{i}}{\mid z_{i}\mid^{2}
(\log\mid z_{i}\mid )^{2}} 
+ \sum_{i=k+1}^{n}\sqrt{-1}dz_{i}\wedge d\bar{z}_{i}
\]
on every compact subset of $\Delta^{n}$ and 
every covariant derivative of $\omega_{M}\mid\Delta^{n}$
 is bounded on every compact subset of $\Delta^{n}$.
\end{definition}

Then by the construction of $\omega_{D}$, we have:

\begin{lemma}(\cite[pp.400-409]{k}) $\omega_{D}$ has
 bounded Poincar\'{e} growth.
on $(X_{\nu},B)$.
\end{lemma}  

Remember the definition of $\varphi$ in Lemma 19. 
Then the following lemma is trivial.

\begin{lemma}  $\varphi^{-1}\!\parallel d\varphi\parallel$
 is uniformly bounded on $X_{\nu} - D$. 
\end{lemma}

We note that $v$ satisfies the following differential inequality.  
\[ 
\Delta_{D}v \geq - tr_{\omega_{D}}\hat{\omega}_{t} \,\,\,
\mbox{on $K\times [0,T)$}, 
\] 
where $\Delta_{D}$ is the Laplacian with respect to $\omega_{D}$. 
Let $h_{D}$ be the solution of the Dirichlet problem 
\[ 
\left\{ \begin{array}{lll} \Delta_{D}h_{D} & = 
-tr_{\omega_{D}}\hat{\omega}_{t} & \mbox{on $K\times [0,T)$} \\ 
h_{D} & = 0 & \mbox{on $\partial K\times [0,T)$} \end{array} \right.  
\] 
Then by the maximum principle, we have 
\[ 
v \leq h_{D} 
 \] 
holds on $K\times [0,T)$.
Hence 
\[ 
\psi_{c} \leq v \leq h_{D}
\]
holds on $K\times [0,T)$ by (12).
Hence by the maximum principle, we have 
\[
\parallel dv\parallel \leq \max \{\parallel d\psi_{c}\parallel 
, \parallel dh_{D}\parallel\} \mbox{on $\partial K\times [0,T)$}.
\]
By the standard boundary estimate for the second order linear
elliptic equations (cf. \cite{g-t}), we see that $h_{D}$ is smooth on $\bar{K}$.  
By using the standard elliptic estimate and Lemma 23, it is 
easy to obtain an estimate for $\parallel dh_{D}\parallel$ on $K$.
But in this case, 
since $\omega_{\varphi} =  \sqrt{-1}\partial\bar{\partial}\varphi$ 
 is a complete K\"{a}hler form of Poincar\'{e} growth, 
we can find a negative constant $b^{\prime}$ 
independent of $c$ and $t$ such that
\[
b^{\prime}\Delta_{D}(\varphi -c) \leq  
-tr_{\omega_{D}}\hat{\omega}_{t}
\]
holds. 
Then by the maximum principle, we see that 
\[
h_{D} \leq b^{\prime}(\varphi - c) \;\;\;\;\mbox{on $K$}
\]
holds. 

Then since $b$ and $b^{\prime}$ are independent of $c$ and $t$,  
by Lemma 25, we have:  

\begin{lemma}
There exists a positive constant $C_{1}^{\prime}$ independent of 
$c \geq c_{0}$ such that
\[
\parallel dv\parallel \leq C_{1}^{\prime}c \;\;\;\; 
\mbox{on $\partial K\times [0,T)$}
\]
where $\parallel\,\,\,\parallel$ is the norm with respect to the 
K\"{a}hler form $\omega_{D}$.  
\end{lemma}

\begin{remark}
 As you  have seen  above, 
in the proof of Lemma 26, the use of $h_{D}$ is not 
unnecessary. 
We can use a $b^{\prime}\psi_{c}$ instead of 
$h_{D}$ from the first.
The reason why we have used $h_{D}$ here is that the method can be  
applicable more general situations. 
\end{remark}

\subsection{$C^{1}$-estimate on $K$}

Let $\pi_{\alpha} : X_{\nu} \longrightarrow {\bf P}^{n}$
 be the generic projection constructed in 6.3.  
And let  
\[
Z =
 \mbox{Re}(\sum_{i}\beta_{i}\frac{\partial}{\partial (z_{i}/z_{0})}),
\]
where $(\beta_{1},\ldots, \beta_{n}) \in {\bf C}^{n}-\{ O\}$.
Then 
\[
\theta = \pi_{\alpha}^{*}(Z)
\]
is a holomorphic differential operator on $X_{\nu}- D$ 
which is meromorphic on $X_{\nu}$.
By operating $\theta$ to (11), we have 
\[
\left\{ \begin{array}{lll}
\frac{\partial (\theta v)}{\partial t} &
 = \tilde{\Delta}(\theta v) -\theta v +
 \theta\log\frac{\hat{\omega}^{n}}{\Omega_{\nu}} & 
\mbox{on $K\times [0,T)$} 
\vspace{5mm} \\
\theta v & =  \theta v &\mbox{on $\partial K\times [0,T)$} \vspace{5mm}\\
\theta v & = 0 & \mbox{on $\partial K\times \{ 0\}$}, 
\end{array}\right.
\]
where $\tilde{\Delta}$ is the Laplacian with repect to 
\[
\omega = \hat{\omega} + \sqrt{-1}\partial\bar{\partial}v.
\]
Then by the maximam principle and Lemma 25, we get 
\[
\parallel\theta v\parallel \leq C_{1}(\theta ,K),
\]
where
\[
C_{1}(\theta ,K) = 
\max \{\sup_{\partial K\times [0,T)}\parallel\theta v\parallel ,
\sup_{K\times [0,T)}
\parallel \theta\log\frac{\hat{\omega}^{n}}{\Omega_{\nu}}\parallel\} .
\]
Then since $\parallel d\log (\hat{\omega}^{n}/\Omega_{\nu})\parallel$ 
 is bounded on $X_{\nu} - D$,  
if we take $m$ sufficiently large and $\pi_{\alpha} (1\leq \alpha \leq m)$ 
properly, we get :
     
\begin{lemma} There exists a positive constant $C_{1}(K)$ 
which depends on  
$c\geq c_{0}$ such that
\[
\parallel dv\parallel \leq C_{1}(K) \;\;\;\; \mbox{on $K\times [0,T)$}
\]
holds.
\end{lemma} 

\vspace{5mm}

The estimate is getting worse if the point goes far from the boundary
 because $\theta$ has a pole along $D$. 

We set
\[
K_{c}(\varepsilon ) = K_{c} - K_{c-\varepsilon}.
\]
Then by the above argument and the construction of $D$ in 6.3, 
we obtain the following estimate.

\begin{lemma} There exists  positive constants $C_{1}$ and $A_{1}$
 independent of $c\geq c_{0}$ such that
\[
\parallel dv\parallel \leq C_{1}c \,\,\,
\mbox{on $K(e^{-A_{1}c})\times [0,T)$}.
\]
\end{lemma}

\begin{remark}
This idea is inspired by the idea in \cite{f}. 
\end{remark}

\subsection{$C^{2}$-estimate on $\partial K$}

In this subsection,  
we follow the argument in \cite[pp. 218-223]{c}
and prove:

\begin{lemma} There exists a positive constant $C_{2}^{\prime}$   
independent of $c\geq c_{0}$ such that 
\[
\parallel \sqrt{-1}\partial\bar{\partial}v\parallel \leq
\exp (C_{2}^{\prime}c) 
\,\,\,
\mbox{on $\partial K\times [0,T)$}
\]
\end{lemma}

\vspace{10mm}

Let $P$ be a point on $\partial K$. Choose a  coordinates 
$z_{1},\ldots ,z_{n}$ with origin at $P$ such that
\begin{enumerate}
\item $dg^{D}_{i\bar{j}}(P) = 0$ and $g^{D}_{i\bar{j}}(P) = \delta_{ij}$.
\item  There exists a positive number $\hat{b}$ such that 
\[
r = \hat{b}(\varphi - c)
\]
satisfies
$r_{z_{\alpha}}(0)= 0$ for $\alpha < n$, $r_{y_{n}}(0) = 0, 
r_{x_{n}} = - 1$, where 
\[
z_{\alpha} = x_{\alpha} + \sqrt{-1}y_{\alpha}
\]
and 
\[
r_{z_{\alpha}} = \frac{\partial r}{\partial z_{\alpha}}
\]
and so on.
\end{enumerate}
We set
$
s_{1}= x_{1},s_{2} = y_{1},\ldots ,s_{2n-3} = x_{n-1}, s_{2n-2} = y_{n-1}, 
s_{2n-1} = y_{n} = s, s^{\prime} = (s_{1},\ldots, s_{2n-1})$. 
By $\partial\bar{\partial}$-Poincar\'{e} lemma, we choose a smooth
function $\phi$ defined on a open neighbourhood $U$ of $P$  such that
\[
\hat{\omega}_{t} = \sqrt{-1}\partial\bar{\partial}\phi 
\]
holds on $U$. 
Let $g$ be a function defined by 
\[
g = \phi + v.
\]

It is clear that to estimate $\sqrt{-1}\partial\bar{\partial}v(P)$ 
is equivalent to $\sqrt{-1}\partial\bar{\partial}g(P)$ because
$\hat{\omega}_{t}$ is uniformly bounded with respect to $\omega_{D}$
on $X_{\nu} - D$ by a constant independent of $t$. . 
 Moreover by Lemma 24, 
 the convariant derivatives of $\hat{\omega}_{t}$ 
 of any order with respect to 
$\omega_{D}$ is uniformly bounded  with respect to
the norm defined by $\omega_{D}$ on $X_{\nu} - D$ 
 by a  constant independent of $t$.
 Then  by Lemma 23  , we may assume that $U$ contains a ball 
of radius $1/2c_{D}$ with center $P$ and any derivatives  of $\phi$ 
of a fixed order with respect to $(z_{1},\ldots ,z_{n})$ 
is bounded by a constant independent of $c\geq c_{0}$,
if we allow $(U,z_{1},\ldots ,z_{n})$ to be a quasi-coordinate. 
Since the estimate is completely local, this does not cause 
any trouble in our estimate in this subsection.  
Hence the $C^{2}$-estimate of $v$ on $\partial K$ is reduced 
completely to the 
$C^{2}$-estimate of $g$ on $\partial K$. 

\begin{sublemma} There exists a positive constant $\dot{C}_{2}$ 
independent of $c\geq c_{0}$ such that
\[
\mid g_{s_{i}s_{j}}(0)\mid \leq \dot{C}_{2} \,\,\,\, i,j \leq 2n -1.
\]
\end{sublemma}
For $r$ near $0$ we may represent $g$ as 
\[
g =  \phi + \sigma r.
\]
Then 
\[
g_{x_{n}}(0) =  \phi_{x_{n}}-\sigma (0),
\]
so that by Lemma 26, $\sigma (0) \leq C_{1}^{\prime}c$. 
Hence 
\[
g_{s_{i}s_{j}}(0) = \phi_{s_{i}s_{j}} + \sigma (0)r_{s_{i}s_{j}}
\]
holds. 
We note that $r_{s_{i}s_{j}} = O(1/c)$ because of the normalization.
Hence we get the sublemma.

\begin{sublemma} 
There exists a positive constant  
$\ddot{C}_{2}$independent of $c\geq c_{0}$  such that
\[
\mid g_{s_{i}x_{n}}(0)\mid \leq \exp(\ddot{C}_{2}c)
\]
holds.
\end{sublemma}

The proof of Sublemma 2 is a little bit technical.

Writing the Taylor expansion of $r$ up to second order we obtain:
\[
r = \mbox{Re}(-z_{n} + \sum a_{ij}z_{i}z_{j}) 
+ \sum b_{i\bar{j}}z_{i}\bar{z}_{j} + O(\mid z\mid^{3}).
\]
Introducing new coordinates of the form 
\[\begin{array}{ll}
z_{n}^{\prime} & = z_{n} - \sum a_{ij}z_{i}z_{j}, \\
z_{k}^{\prime} &= z_{k} \,\,\mbox{for $k\leq n-1$},
\end{array}
\]
we can write
\begin{equation}
r = -\mbox{Re}(z^{\prime}_{n}) +
 \sum c_{i\bar{j}}z_{i}^{\prime}\bar{z}_{j}^{\prime} + O(\mid z\mid^{3}).
\end{equation}
It is clear that $(c_{i\bar{j}})$ is positive definite. 

We define $T_{i}$ in a neighbourhood of $0$ by
\[
T_{i} = \frac{\partial}{\partial s_{i}} - \frac{r_{s_{i}}}{r_{x_{n}}}
\frac{\partial}{\partial x_{n}}, \mbox{for $i = 1,\ldots ,2n - 1$};
\]
then $T_{i}r = 0 $ nad we have $T_{i}(g - \phi ) = 0$ on 
$r = 0$. 

We show that for suitable $\varepsilon > 0$, in the region 
\[
S_{\varepsilon} = \{ x\in U\mid r(x)\leq 0, x_{n}\leq \varepsilon\},
\]
where $U$ is a neighbourhood of the origin, we set
\[
w = \pm T_{i}(g -\phi) + (g_{s} -\phi_{s})^{2} - Ax_{n} + B\mid z\mid^{2} 
\]
We claim :

\vspace{5mm}

(a) For $B$ sufficiently large, $\tilde{L}w \geq 0$;

(b)On $\partial S_{\varepsilon}$,
 if $A$ is sufficiently large,
$w\leq 0$ holds.

\vspace{5mm}

To prove  (a) set  
\[
a = - r_{s_{i}}/r_{x_{n}}
\]
and consider (we use summation convention)
\begin{equation}
\tilde{L}(T_{i}g) = T_{i}\log\Psi (z,g(z)) + 
g^{p\bar{q}}a_{p}g_{x_{n}\bar{q}} + g^{p\bar{q}}a_{\bar{q}}g_{x_{n},p}
 + g^{p\bar{q}}a_{p\bar{q}}g_{x_{n}},
\end{equation}
where
\[
\tilde{L} = \tilde{\Delta} -\frac{\partial}{\partial t},
\]
and 
\[
\Psi (z,g(z)) =
 \frac{\Omega_{\nu}}{(\sqrt{-1})^{n}dz_{1}\wedge d\bar{z}_{1}\wedge
\ldots\wedge dz_{n}\wedge d\bar{z}_{n}} \exp (-g\phi ).
\]
Obeserve that $g^{p\bar{q}}g_{n\bar{q}} = \delta^{p}_{n}$ and that
\[
\frac{\partial}{\partial x_{n}} = 2\frac{\partial}{\partial z_{n}} + 
\sqrt{-1}\frac{\partial}{\partial s}
\]
so that
\[
g_{x_{n}\bar{q}} = 2g_{n\bar{q}} + \sqrt{-1}g_{s\bar{q}}.
\]
Thus  the second term on the right-hand side of (16) is of the 
form 
\[
a_{n} + g^{p\bar{q}}a_{p}g_{t\bar{q}} = 
O(1+ (\sum g^{i\bar{i}})^{1/2}
(g^{p\bar{q}}g_{ps}g_{\bar{q}s})^{1/2}).
\]
A similar estimate holds for the third term on the right of (16)
 while the forth term is $O(\sum g^{i\bar{i}})$.
Thus by Lemma 19 and the arithmetic-geometric mean inequality, we
have
\[
\pm \tilde{L}T_{i}g \leq - C\Psi^{-1/n} - g^{p\bar{q}}g_{ps}g_{\bar{q}s}.
\]
Further
\[ \begin{array}{ll}
\tilde{L}(g_{s}-\phi_{s})^{2} & = 2g^{p\bar{q}}g_{ps}g_{\bar{q}s} +
2(g_{s} - \phi_{s})(\partial_{s}\log \Psi - \tilde{L}\phi_{s}) \\
 & \geq 2g^{p\bar{q}}g_{ps}g_{\bar{q}s} - C\Psi^{-1/n}
\end{array}
\]
holds on $S_{\varepsilon}$ by  the $C^{1}$-estimate (Lemma 28) 
 and the arithmetic-geometric mean
 inequality (if we take $\varepsilon$ sufficiently small). 
Hence we find
\[ \begin{array}{ll}
\tilde{L}w & \geq B\sum g^{i\bar{i}} - C\Psi^{-1/n} \\
 & \geq 0 \,\,\,\,\,\mbox{on $S_{\varepsilon}\times [0,T)$} ,
\end{array}
\]
if $B$ is sufficiently large. 
$B$ depends on the $C^{1}$-estimate of $v$ 
 on $S_{\varepsilon}\times [0,T)$ and which is uniform
 with the weight $\varphi$ by Lemma 27.
This completes the proof of (a).
 
To prove (b), consider first $\partial S_{\varepsilon}\cap \partial K$.  
Here we write $x_{n} = \rho (s_{1},\ldots ,s_{2n-1})$ and
from (13), we deduce that 
\begin{equation}
\rho (t) = \sum_{ij < 2n}b_{i\bar{j}}s_{i}s_{j} + O(\mid s^{\prime}\mid^{3})
\end{equation}
with $(b_{i\bar{j}})$ positive definite.   
Thus on $\partial K$ near $O$ we have 
\[
x_{n}\geq a\mid z\mid^{2}, 
\]
where $a$ is uniformly bounded from below by a positive constant 
 times $1/c$ where by the construction of $K$. 
Also,
\[
g(s,\rho (s)) = \phi (s,\rho (s)),
\]
so that 
\[
\mid g_{s} - \phi_{s}\mid^{2} \leq C\mid s\mid^{2} \leq C\rho .
\]
Taking $A$ large we obtain (b).
By Lemma 28. if we take a positive constant $\ddot{C}_{2}$ 
sufficiently
large,    
we  may assume that $A$ is bounded from above by a constant times 
$\exp (\ddot{C}_{2}c)$ for some positive constant $\ddot{C}_{2}$  
 independent of $c \geq c_{0}$.
By the maximum principle and (a),(b), 
\[
w\leq 0 \;\;\;\;\;\;\mbox{on $S_{\varepsilon}$} 
\]
holds. 
 
In view of the  maximum principle, we have  
\[ 
\mid (T_{i}g)_{x_{n}}(0)\mid \leq A. 
\] 
This completes the proof of Sublemma 2. \vspace{5mm} 

Using, still the special coordinate above we see that
 to finish our proof of Lemma 29 , we have only to establish the estimate
\begin{equation}
\mid g_{x_{n}x_{n}}(O)\mid \leq  \exp (\acute{C}_{2}c)
\end{equation}
for some constant $\acute{C}_{2}$ independent of $c \geq c_{0}$. 
By the previous estimates: 
\begin{eqnarray*}
\mid g_{s_{i}s_{j}}(O)\mid \leq \dot{C}_{2} & (1\leq i,j\leq 2n -1),\\
 \mid g_{s_{i}x_{n}}\mid 
\leq \exp (\ddot{C}_{2}c) &  (1\leq i\leq 2n- 1),
\end{eqnarray*}
 it suffices to prove
\[
\mid g_{n\bar{n}}(O)\mid \leq \exp (C_{2}^{\prime}c)
\]
for some $\hat{C}_{2}$ independent of $c\geq c_{0}$. 
We may solve the equation 
\[
\det (g_{i\bar{j}})(O) = \frac{\Omega_{\nu}}{\omega_{D}^{n}}(O)
\]
for $g_{n\bar{n}}(O)$.   
 Then  since there exists a positive  constant $C^{\flat}$ such that 
\[
\exp (-C^{\flat}c) \leq \frac{\omega_{D}^{n}}{\Omega_{\nu}} 
\leq \exp (C^{\flat}c) 
\]
holds on $X_{\nu} - D$, 
 we see that (16) follows from (17) provided we know
the following sublemma.  
\begin{sublemma}(\cite[pp. 221-223]{c}) 
The $(n -1)$ by $(n -1)$ matrix
\[
(g_{z_{\alpha}\bar{z}_{\beta}}(O))_{\alpha ,\beta < n} \geq 
C_{3}(\frac{\Omega_{\nu}}{\omega_{D}^{n}})^{\frac{1}{n}}\frac{1}{c}I
\]
for some $C_{3}$; here $I$ is the $(n -1)$ by $(n -1)$ identity matrix.  
$C_{3}$ is independent of $c\geq c_{0}$ 
\end{sublemma} 

The proof of this sublemma is very technical. 

After subtraction of a linear function we may assume that 
$\phi_{s_{j}}(O) = 0, j\leq 2n -1$. 
To prove Sublemma 3, it suffices to prove 
\[ 
\sum_{\alpha ,\beta <n}
\gamma_{\alpha}\bar{\gamma}_{\beta}g_{z_{\alpha}\bar{z}_{\beta}}(O)\geq C_{3}\mid\gamma\mid^{2}
\]
which we shall do for $\gamma = (1,0,\ldots ,0)$. 
We shall show that
\begin{equation}
 g_{1\bar{1}}(O)\geq C_{4},
\end{equation}
where $C_{4}$ is a positive constant. 
Let $\tilde{g} = g - \lambda x_{n}$ with $\lambda$ so chosen that
\[
(\frac{\partial^{2}}{\partial s_{1}^{2}}+
\frac{\partial^{2}}{\partial s_{2}^{2}})\tilde{g}(s_{1},\ldots ,s_{2n-1},
\rho (s_{1},\ldots ,s_{2n-1})) = 0 \;\;\; \mbox{at $O$},
\]
i.e.
\begin{equation}
 0 = g_{1\bar{1}}(O)+\tilde{g}_{x_{n}}(O)\rho_{1\bar{1}}(O) =
 g_{1\bar{1}}(O) + (g_{x_{n}}(O) - \lambda )\rho_{1\bar{1}}(O).
\end{equation}
Using the fact that any real homogeneous cubic polynomial 
in $(s_{1},s_{2})$ admits the unique decomposition
\[
\mbox{Re}(\alpha (s_{1}+\sqrt{-1}s_{2})^{3}+\beta 
(s_{1}+\sqrt{-1}s_{2})(s_{1} +\sqrt{-1}s_{2})^{2}),
\]
we find on expanding $\tilde{g}\mid_{\partial K\cap U}$ 
in a Taylor series, in $s_{1},\ldots ,s_{2n-1}$,
\[ \begin{array}{ll}
\tilde{g}\mid_{\partial K\cap U} =
  \mbox{Re}\sum_{2}^{n-1}a_{j}z_{1}\bar{z}_{j} 
+\mbox{Re}(az_{1}s) & + \mbox{Re}(p(z_{1},\ldots ,z_{n-1}) 
+ \beta z_{1}\mid z_{1}\mid^{2})  \\
 & + O(s_{3}^{2}+ \ldots + s_{2n-1}^{2}), 
\end{array}
\]
where $p$ is a holomorphic cubic polynomial.

With the aid of (13), we may replace the term $\beta z_{1}\mid z_{1}\mid^{2}$ 
 to $(\rho_{z_{1}\bar{z}_{1}}(O))^{-1}\beta z_{1}x_{n}$,
if we change the $a_{j},a$ and $p$. 
Thus by changing the $a_{j},a$ and $p$ appropriately 
we may obtain the inequality:
\[
\tilde{g}\mid_{\partial K\cap U}\leq \mbox{Re}\,\,p(z) + 
\mbox{Re}\sum_{2}^{n}a_{j}z_{1}\bar{z}_{j} + 
C\sum_{j=2}^{n}\mid z_{j}\mid^{2}.
\]
Let $\check{g} = \tilde{g} - \mbox{Re}\,p(z)$ and observe that
$\check{g}$ satisfies 
\[
\det (\check{g}_{j\bar{k}})
 = \det (\tilde{g}_{j\bar{k}}) = \Psi (z,g(z)).
\]
Recall that $\Psi (z,g(z)) \geq \delta > 0$  on a neighbourhood of 
$O$, where $\delta$ depends on $\Omega_{\nu}/\omega_{D}^{n}(O)$
 and the $C^{1}$-estimate of $v$  on the neighbourhood. 
With $\varepsilon$ small we see that in the set $S_{\varepsilon}$, we
have $\Psi (z,g(z)) \geq \delta$.   
Let 
\[
h = -\delta_{0}x_{n} + \delta_{1}\mid z\mid^{2} 
+\frac{1}{B}\sum^{n}_{2}\mid a_{j}z_{1} + Bz_{j}\mid^{2}.
\]
We wish to show that with the suitable choice of 
$\delta_{0},\delta_{1}, B > 0$ we have $h\geq \check{g}$
 on $\partial S_{\varepsilon}$. 
First observe that if $B$ is sufficiently large and $\delta_{0}$
 so small that $-\delta_{0}x_{n} +\delta_{1}\mid z\mid^{2}\geq 0$ 
on $\partial S_{\varepsilon}\cap\partial K$ 
 (the dependence of $\delta_{0}$ and $\delta_{1}$ is controlled
 by the Levi form of $\partial K$).
By Lemma 28, if we take $B$ sufficiently large, we have
\[
\check{g} \leq h \;\;\;\;\;\mbox{on $\partial S_{\varepsilon}$.}
\]  
The function $h$ is plurisubharmonic and the lowest eigenvalues of 
the complex Hessian $(h_{i\bar{j}})$ are bounded independently by
 $\delta_{1}$ while the other eigenvalues are bounded 
independently of $\delta_{1}$.    

Hence choosing $\delta_{1}$ equal to small const. times $\delta^{1/n}$
\[
\det (h_{i\bar{j}}) \leq \delta  \,\,\,\, \mbox{in $S_{\varepsilon}$}
\]
holds. 
By the mximum principle 
\[
\check{g} \leq h \;\;\;\;\;\;\mbox{on $S_{\varepsilon}$}
\]
holds. 
Hence by the maximum principle
\[
\check{g}_{x_{n}}(O)\leq h_{x_{n}}(O) = - \delta_{0}. 
\]
The desired inequality follows from (19). 
\vspace{5mm}
This completes the proof of Lemma 29.  

\subsection{$C^{2}$-estimate on $K$}  

Using the $C^{2}$-estimate on $\partial K$, 
we shall obtain a $C^{2}$-estimate inside $K$. 
The method here is the same as in \cite{t}.

Let $H$ be a smooth function on $X_{\nu} - D$ defined by
\[
H = \exp (\delta_{\nu})
(\prod_{k}\parallel\lambda_{k}\parallel
\prod_{k}(\log
\frac{1}{\parallel\lambda_{k}\parallel})^{-1})^{\varepsilon}
\]
where $\delta_{\nu}$ is the one in 4.2 and
 $\parallel\lambda_{k}\parallel$'s  are  the ones in 4.6  
and $\varepsilon$ is a sufficiently small positive
 number such that
\[
\omega_{H} = \hat{\omega}_{t} + \sqrt{-1}\partial\bar{\partial}\log H 
\]
is a complete K\"{a}hler form on $(X_{\nu} - D)\times [0,\infty]$ 
 which is  quasi-isometric to $\omega_{D}$ on $X_{\nu} - D$,
 i.e., there exists a positive constant $C(D,H) > 1$ such that
\[
\frac{1}{C(D,H)}\omega_{H}\leq \omega_{D}\leq C(D,H)\omega_{H}
\]
holds on $X_{\nu} - D$.
We note that $\omega_{H}$ have bounded Poincar\'{e} growth so that
the bisectional curvature of $\omega_{H}$ is bounded between 
two constants uniformly on $(X_{\nu}- D)\times [0,\infty ]$. 

We set 
\[
v_{H} = v - \log H = u - (1 -e^{-t^{4}})\xi_{\nu} - \log H, 
\]
\[
\Omega_{H} =  H\cdot\Omega_{\nu}.
\]
Then $v_{H}$ satisfies the equation
\[
\frac{\partial v_{H}}{\partial t} = 
\log\frac{(\omega_{H}+
\sqrt{-1}\partial\bar{\partial}v_{H})^{n}}{\Omega_{H}}-v_{H} \;\;\; 
\mbox{on $K\times[0,T)$}.
\]
By Lemma 20 and Lemma 22 $v_{H}$ satisfies the $C^{0}$-estimate:

\begin{lemma}
For every sufficiently small positive number $\varepsilon$ 
\[\begin{array}{ll}
v_{H}\geq &  C_{0}^{-}(\varepsilon) -\log H  + 
(1- e^{-t^{4}})\varepsilon\delta_{\nu}  \\
v_{H}\leq & C_{0}^{+}(1-e^{-t}) - (1 -e^{-t^{4}})\xi_{\nu} -\log H 
\end{array}
\]
holds on $K\times [0,T)$, where $C_{0}^{+},C_{0}^{-}(\varepsilon)$ 
are constants in Lemma 20 and 22 respectively.
\end{lemma} 

We have the following lemma. 

\begin{lemma}(\cite[Lemma 3.2]{t})
\begin{eqnarray*}
H^{-C}e^{Cv}(\tilde{\Delta}-\frac{\partial}{\partial t})
(e^{-Cv}H^{C}tr_{\omega_{H}}\omega )  \geq 
\;\;\;\;\;\;\;\;\;\;\;\;\;\;\;\;\;\;\;\;\;\;\;\;\;\;\;\;\;\;\;\;\;\;\;\\
(-\Delta_{H}\log\frac{\omega_{H}^{n}}{\Omega_{H}} 
- n^{2}\inf_{i\neq j}R_{i\bar{i}j\bar{j}}^{H}-n) 
+ C(n -\frac{1}{C}-\frac{\partial v_{H}}{\partial t})
tr_{\omega_{H}}\omega - \;\;\;\;\;\;\;\;\;\;\;\;\;\; \\
e^{-t}tr_{\omega_{H}}((\omega_{0}-\omega_{\infty})\cdot\omega )  
+ (C+\inf_{i\neq j}R_{i\bar{i}j\bar{j}}^{H})
\exp(\frac{1}{n-1}(-\frac{\partial v_{H}}{\partial t}
-v_{H}+\log\frac{\omega_{H}^{n}}{\Omega}))
(tr_{\omega_{H}}\omega )^{\frac{n}{n-1}} ,
\end{eqnarray*}
holds on $K\times [0,T)$,  where 
$tr_{\omega_{H}}((\omega_{0}-\omega_{\infty})\cdot\omega)$
 is defined as in Lemma 3.4,  
$\inf_{i\neq j}R^{H}_{i\bar{i}j\bar{j}}$ denotes the infimum
 of the bisectional curvarue of $\omega_{H}$ on 
$(X_{\nu} - D)\times [0,\infty ]$ and
$C$ is a positive constant such that 
\[
C + \inf_{i\neq j}R^{H}_{i\bar{i}j\bar{j}} > 1 
\]
holds.
\end{lemma} 

\vspace{10mm}

The proof of this lemma is the same as one of  Lemma 3.2 in \cite{t}.
Hence we omit it.

\begin{lemma} 
If we take $C$ sufficiently large, then 
there exists a positive constant $C_{2}$ independent of 
$c\geq c_{0}$ such that 
\[
H^{C}tr_{\omega_{D}}\omega \leq C_{2} \,\,\,\mbox{on $K\times [0,T)$}
\]
holds. 
\end{lemma}
{\em Proof}. 
By the definition $H$ has zero  of order  
at least $r_{j}^{\nu}/2$ along $R_{j}^{\nu}$ (cf. 6.2). 
 Then by Lemma 29, if we take $C$ sufficiently large, there exists 
a constant $\tilde{C}_{2}$ independent of $c\geq c_{0}$  such that 
\[
H^{C}tr_{\omega_{D}}\omega \leq \tilde{C}_{2} \,\,\,
\mbox{on $\partial K\times [0,T)$}.
\]
holds.
Suppose  $e^{-Cv}H^{C}tr_{\omega_{H}}\omega$ takes its maximum at
$P_{0} \in K\times \{ t_{0}\} (t_{0}\in [0,T))$ then by Lemma 31,
we have
\[
(tr_{\omega_{H}}\omega )(P_{0}) \leq \hat{C}_{2}
\]
for a positive constant $\hat{C}_{2}$ independent of $c\geq c_{0}$
 and $C$(if it is sufficiently large).
Hence in this case we have
\[
(H^{C}tr_{\omega_{H}}\omega )(P) 
\leq H^{C}(P_{0})\exp(-C(v(P_{0})-v(P)))\hat{C}_{2} \;\;\;\;
\mbox{on $K\times [0,T)$} 
\]
holds. 
By Lemma 30 
(since in Lemma 30, we can take $\varepsilon$ arbitrarily small),
 $H\exp (-v) = \exp (-v_{H})$ is uniformly bounded from above on
 $X_{\nu} - D$. 
Hence if we change $\hat{C}_{2}$, if necessary,
by the maximum principle for parabolic equations, we may assume that
\[
(H^{C}tr_{\omega_{H}}\omega )(P)\leq \exp (Cv(P))\hat{C}_{2}
\]
holds on $K\times [0,T)$.  
We note that $\omega_{D}$ and $\omega_{H}$ are quasi-isometric 
on $X_{\nu} - D$. 
Hence by Lemma 22, if we take $C$ sufficiently large, 
there exists a positive constant $C_{2}$ independent of
 $c\geq c_{0}$ such that    
\[
H^{C}tr_{\omega_{D}}\omega \leq C_{2} \,\,\,\mbox{on $K\times [0,T)$}
\] 
holds.       
Q.E.D. 
   
\vspace{10mm}

By \cite{tr}, the higher order interior estimate 
on $K\times [0,T)$ follows.  
As for the boundary estimate of $u$, we just need to follow the 
argument in $\cite{c}$.   This completes the proof of Theorem 10. 

\subsection{Costruction of $\omega_{E}$}

Let $u_{\infty}$ be as in Theorem 10.  
Then by the construction
\[
\omega_{K}  = \omega_{\infty} + \sqrt{-1}\partial\bar{\partial}u_{\infty}
\]
is a K\"{a}hler-Einstein form on $\bar{K} = \bar{K}_{c}$. 
We may assume that $c > 1$.  
Let us take the exahustion $\{ K_{lc}\}_{l=1}^{\infty}$
 of $X_{\nu}- D$ and let
\[
\omega_{l}^{\nu} = \omega_{K_{lc}}.
\]
Then by  Lemma 20 and Lemma 32  and the regularity theorem 
in \cite{tr}, we have the following lemma.

\begin{lemma}  There exists a subsequence of 
$\{ \omega_{l}^{\nu}\}_{l=1}^{\infty}$
 which converges uniformly on every compact subset of $X_{\nu}- D$
 in $C^{\infty}$-topology to a K\"{a}hler-Einstein form
 $\omega^{\nu}$ on $X_{\nu} - D$.
\end{lemma}

\vspace{10mm}

Although $\omega^{\nu}$ is a K\"{a}hler-Einstein form,
 it is not enough good for our purpose. 

Let us consider the linear system  
$\mid m!\nu L\mid$ and construct $\xi_{m!\nu}$ as before.  
We  denote $\xi_{m!\nu}$ by $\xi^{(m)}$ for simplicity. 
Since $X_{m!\nu} - F^{m!\nu}$  are all biholomorphic 
to $X - \mbox{SBs}(L)$, we may consider
 $\{\xi^{(m)}\}$ as a family of functions on 
$X_{\nu} - D$. 
Let us denote $X_{m!\nu}$ by $X^{(m)}$ 
for simplicity and define $X_{(1)} = X$. 

\begin{lemma} We can construct $\{\xi^{(m)}\}$ so that
\[
\xi^{(1)}\leq \xi^{(2)}\leq\ldots\leq \xi^{(m)} \leq \xi^{(m+1)}\leq\ldots 
\]
holds on $X_{\nu} - D$.
\end{lemma}
{\em Proof}. Let ${\cal I}_{\mu}$ denote the ideal sheaf of
 the base scheme $\mbox{Bs}\mid\mu L\mid$. Then clearly
\[
{\cal I}_{\mu_{1}+\mu_{2}}\hookrightarrow 
{\cal I}_{\mu_{1}}\otimes_{{\cal O}_{X}}{\cal I}_{\mu_{2}} 
\]
holds. 
Hence inductively we can construct a morphism 
\[
\mu_{m} : X_{(m)}\longrightarrow X_{(m-1)} \,\,\,(m\geq 2)
\]
such that
\begin{enumerate}
\item $f_{(m)} := \mu_{m}\circ\cdots\circ\mu_{1}: X_{(m)}
\longrightarrow X_{(0)}(= X)$ is a resolution of \\
 $\mbox{Bs}\mid m!\nu L\mid$. 
\item The fixed part of $\mid f_{(m)}^{*}(m!L)\mid$ 
is a divisor with normal crossings on $X_{(m)}$. 
\end{enumerate}

 An explicit construction of $\{\xi^{(m)}\}_{m=1}^{\infty}$
 is as follows. 
Let $V^{(1)} = \{\eta_{i}^{(1)}\}_{i=0}^{N(1)}$ be a basis of 
$H^{0}(X,{\cal O}_{X}(\nu L))$.
By induction for each $m\geq 1$, we can construct a finite  subset 
\[
V^{(m)} = \{ \eta^{(m)}_{1},\ldots ,\eta^{(m)}_{N(m)}\} 
\]
in $H^{0}(X,{\cal O}_{X}(m!\nu L)$ with the following
properties. 
\begin{enumerate}
\item $V^{(m)}$ spans $H^{0}(X,{\cal O}_{X}(m!\nu L))$.
\item $V^{(m)}$ contains all the elements of the form:
\begin{eqnarray*}
(\eta^{(m-1)}_{i_{1}})^{\otimes a_{1}}\otimes\cdots\otimes 
(\eta^{(m-1)}_{i_{m}})^{\otimes a_{m}}, \\ 
\sum_{i=1}^{m}a_{i} = m, a_{i}\geq 0,      
0\leq i_{1}<\ldots < i_{m}\leq N(m-1). 
\end{eqnarray*}
\end{enumerate} 
Now we set 
\[
\xi^{(m)} = \frac{1}{m!\nu}\log (\sum^{N(m)}_{i=0}
\frac{(\sqrt{-1})^{m!\nu n^{2}}\eta^{(m)}_{i}\wedge\bar{\eta}^{(m)}_{i}}
{\Omega^{m!\nu}}).   
\]
We may consider $\xi^{(m)}$ as a function on $X_{\nu} - D$.
Then by the construction 
\[
\xi^{(1)}\leq \xi^{(2)}\leq\ldots\leq\xi^{(m)}\leq\ldots
\]
holds and 
\[
\omega_{\infty} + \sqrt{-1}\partial\bar{\partial}\xi^{(m)} 
\]
is a smooth semipositive form on $X_{(m)}$ and positive on 
$X_{\nu} - D$. 
This completes the proof of the lemma. 
Q.E.D. 

\vspace{10mm}

Now we consider the following Dirichlet problem.
\[
\left\{ \begin{array}{lll}
\frac{\partial u^{(m)}}{\partial t} & = 
\log\frac{(\omega_{t}+
\sqrt{-1}\partial\bar{\partial}u^{(m)})^{n}}{\Omega} - u^{(m)} &
\mbox{on $K\times [0,T_{m})$} \vspace{5mm}\\
u^{(m)} & = (1- e^{-t^{4}})\xi^{(m)}  
& \mbox{on $\partial K\times [0,T_{m})$}\vspace{5mm} \\
u^{(m)} & = 0 & \mbox{on $K\times \{ 0\}$,}
\end{array}
\right. 
\]
where $T_{m}$ is the maximal existence time 
for smooth solution on $\bar{K}$. 

\begin{lemma}The followings are true.
\begin{enumerate}
\item  $T_{m}$ is infinite and 
$u^{(m)}_{\infty} = \lim_{t\rightarrow\infty}u^{(m)}$ exists 
in $C^{\infty}$-topology on $\bar{K}$. 
\item $\omega^{(m)}_{1}:= 
\omega_{\infty} + \sqrt{-1}\partial\bar{\partial}u^{(m)}_{\infty}$
 is a K\"{a}hler form on $K$.
\item If we define a sequence of K\"{a}hler forms
 $\{\omega^{(m)}_{l}\}_{l=1}^{\infty}$ in the same manner 
as the definition of  $\{\omega^{\nu}_{l}\}_{l=1}^{\infty}$ above,
 then there exists a subsequence of $\{\omega_{l}^{(m)}\}$
 which converges uniformly on every compact subset of 
$X_{\nu}- D$ uniformly in $C^{\infty}$-topology. 
\end{enumerate} 
\end{lemma}
{\em Proof}. 
The only difference between the above equation and the equation (10) is
that $\nabla^{k}\xi^{(m)} (k\geq 1)$ is  bounded 
with respect to $\omega_{D}$ with weights different from before.  
It is easy to find such weights. In fact, there exists a positive
 constant $C(m,k)$ depending only on $m$ and $k$ such that
$H^{C(m,k)}\parallel\nabla^{k}\xi^{(m)}\parallel$ is bounded
 by a positive constant on $X_{\nu}- D$.  
Hence the previous argument is valid with 
some weight with respect to $H$.
  Hence the first assertion is trivial. 

Then by replacing $K$ to  $K_{lc}$, we get a sequence of
 K\"{a}hler-Einstein form $\{ \omega_{l}^{(m)}\}_{l=1}^{\infty}$
 which are defined on $K_{lc}$ respectively
 ($\omega^{(1)}_{l} = \omega^{\nu}_{l}$).  

We would like to find a subsequence of 
$\{ \omega_{l}^{(m)}\}_{l=1}^{\infty}$ which converges
 in $C^{\infty}$-topology on every compact subset of $X_{\nu} - D$.

For the first,  replacing $\xi_{\nu}$ by $\xi^{(m)}$, 
completely analogous estimate as Lemma 22 holds for $u^{(m)}$ with  
the perturbation $\delta^{(m)}$ completely analogous to $\delta_{\nu}$.  
 As for the $C^{2}$-estimate, by the proof of Lemma 35, Lemma  32
 also holds for $\omega_{l}^{m}$ if we replace $C$ and $C_{2}$
 to appropriate  constants independent of $l$.
The rest of the proof is the same as in the proof of Lemma 33. 
Q.E.D. 

\vspace{10mm}

Taking subsequence, if necessary, 
we obtain 
$\omega^{(m)}$ as before, where $\omega^{(0)} = \omega^{\nu}$.  
        
We shall consider  $\omega^{(m)}$ as  
a $d$-closed positive $(1,1)$-current on $X$  by 
\[
\omega^{(m)} := \omega_{\infty} +
 \sqrt{-1}\partial\bar{\partial}\log\frac{(\omega^{(m)})^{n}}{\Omega},
\]
where $\partial\bar{\partial}$ is taken in the sense of a current.
This definition is well defined by the $C^{0}$-estimate (Lemma 22)
 and clearly 
\[
[\omega^{(m)}] = 2\pi c_{1}(L)
\]
holds. 
Now we want to show that
\begin{proposition}
\[
\omega_{E} = \lim_{m\rightarrow\infty}\omega^{(m)}
\]
exists 
in the sense of a  $d$-closed positive $(1,1)$-current on $X_{\nu}$. 
\end{proposition}

\begin{lemma} 
There exists a positive constant $C_{3}$ independent of $m$ such that 
\[
(\omega^{(m)})^{n}\leq  C_{3}\omega_{D}^{n} \,\,\,\,\,\mbox{on $X_{\nu} - D$}
\]
holds.
\end{lemma}
{\em Proof}. 
 Since $\omega_{D}$ is a {\em complete} K\"{a}hler-Einstein form on
 $X_{\nu}- D$, by applying the maximum principle to the function 
 \[
 \log\frac{(\omega^{(m)})^{n}}{\omega_{D}^{n}},
 \]
we obtain the lemma.  Q.E.D.

\vspace{10mm}

Hence $\{ (\omega^{(m)})^{n}\}$ is uniformly bounded from above. 
For the next we shall show:

\begin{lemma} For every $m\geq 1$, we have that 
\[
(\omega^{(m)})^{n} \leq (\omega^{(m+1)})^{n}  \,\,\,\,\,
\mbox{on $X_{\nu} - D$}
\]
holds.
\end{lemma} {\em Proof}. $u^{(m+1)} - u^{(m)}$ satisfies 
the following equations: 
\[ \left\{ \begin{array}{lll} 
\frac{\partial (u^{(m+1)}-u^{(m)})}{\partial t}
 & =
\log \frac{(\omega_{t}+\sqrt{-1}\partial\bar{\partial}u^{(m+1)})^{n}}
{(\omega_{t} +\sqrt{-1}\partial\bar{\partial}u^{(m)})^{n}} 
- (u^{(m+1)}-u^{(m)}) &\mbox{on $K\times [0,\infty )$} \vspace{5mm}\\
u^{(m+1)}- u^{(m)}& = (1 - e^{-t^{4}})(\xi^{(m+1)}-\xi^{(m)}) &
\mbox{on $\partial K\times [0,\infty )$} \vspace{5mm} \\ 
u^{(m+1)}- u^{(m)} & = 0 &\mbox{on $K\times\{ 0\}$.}
 \end{array}\right.
\] 
Since
\[
\log\frac{(\omega_{t}+ \sqrt{-1}\partial\bar{\partial}u^{(m+1)})^{n}}
{(\omega_{t}+\sqrt{-1}\partial\bar{\partial}u^{(m)})^{n}} = 
\int_{0}^{1}\Delta_{a}^{(m,m+1)}(u^{(m+1)}-u^{(m)})da, 
\]
where $\Delta^{(m,m+1)}_{a} (a\in [0,1])$ is
 the Laplacian with respect to the 
K\"{a}hler form
\[
\omega_{t} + \sqrt{-1}\partial\bar{\partial}
\{ (1-a)u^{(m)}+au^{(m+1)}\} ,
\]
we see that this equation is of parabolic type. 
We note that $\xi^{(m+1)}\geq \xi^{(m)}$ on $X_{\nu} - D$
 by the construction.
Then by the maximum principle we obtain 
\[
u^{(m)} \leq u^{(m+1)}  \mbox{on $K\times [0,\infty )$}.
\]                       
We set
\[
u^{(m)}_{\infty} = \lim_{t\rightarrow\infty}u^{(m)}.
\]

Then we have 
\[
\omega^{(m)}_{1}  
 =  \omega_{\infty} + \sqrt{-1}\partial\bar{\partial}u_{\infty}^{(m)}
\] 
and
\[
(\omega_{1}^{(m)})^{n} = \exp (u^{(m)}_{\infty})\Omega
\]
on $K$. 
Hence we see that
\[
(\omega_{1}^{(m)})^{n}\leq (\omega_{1}^{(m+1)})^{n}
\]
holds on $K$.
By replacing $K$ to $K_{lc}$ and repeating the same argument, we 
see that
\[
(\omega_{l}^{(m)})^{n}\leq (\omega_{l}^{(m+1)})^{n}
\]
holds on $K_{lc}$. 
By letting $l$ tend to infinifty, we completes the proof of the
lemma.  Q.E.D.

\vspace{10mm}

Hence $\{ (\omega^{(m)})^{n}\}_{m=1}^{\infty}$ is monotone increasing 
and bounded from above uniformly on every compact subset of
 $X_{\nu} - D$ 
\[
\omega_{E}^{n}:=\lim_{m\rightarrow\infty}(\omega^{(m)})^{n}
\]
exists. 

We shall define a $d$-closed positive $(1,1)$ current 
$\omega_{E}$ on $X$ by 
\[
\omega_{E} = 
 \sqrt{-1}\partial\bar{\partial}\log\,\omega_{E}^{n}.
\]
$\omega_{E}$ is well defined by the $C^{0}$-estimates
 in the last subsection.
Then by the construction it is clear that $[\omega_{E}] = 2\pi c_{1}(K_{X})$ and $(K_{X},(\omega_{E}^{n})^{-1})$  is an AZD. 
This completes the proof of Theorem 7.

\end{document}